# Random Vector Representation of Continuous Functions and Its Applications in Quantum Mechanics

Hong-Xing Li  Wei Zhou  Hong-Hai Mi

Research Center of Applied Mathematics and Interdisciplinary Sciences
Beijing Normal University, Zhuhai, 519087, China
Emails: lihx@dlut.edu.cn, 9393560@qq.com, mhhygx@126.com

**Abstract** The relation between continuous functions and random vectors is revealed in the paper that the main meaning is described as, for any given continuous function, there must be a sequence of probability spaces and a sequence of random vectors where every random vector is defined on one of these probability spaces, such that the sequence of conditional mathematical expectations formed by the random vectors uniformly converges to the continuous function. This is random vector representation of continuous functions, which is regarded as a bridge to be set up between real function theory and probability theory. By means of this conclusion, an interesting result about function approximation theory can be got. The random vectors representation of continuous functionsis is of important applications in physics. Based on the conclusion, if a large proportion of certainty phenomena can be described by continuous functions and random phenomenon can also be described by random variables or vectors, then any certainty phenomenon must be the limit state of a sequence of random phenomena. And then, in the approximation from a sequence of random vectors to a continuous function, the base functions are appropriately selected by us, an important conclusion for quantum mechanics is deduced: classical mechanics and quantum mechanics are unified. Particularly, an interesting and very important conclusion is introduced as the fact that the mass point motion of a macroscopical object possesses a kind of wave characteristic curve, which is called wave-mass-point duality.

## 1. Introduction

From a physical point of view, continuous functions can describe a large proportion of certainty phenomena. For example, the trajectory of a projectile motion can be expressed as a continuous function. However random vectors can describe a lot of random phenomena in natural world. So if we consider the connection between some certainty phenomena and some random phenomena, we should research the relation between continuous functions and random vectors.

In the paper, we obtain a conclusion: for any given continuous function $f(x) \in C[a,b]$, there must be a sequence of probability spaces $\{(\Omega, \mathscr{F}, P_n)\}$ and a sequence of random vectors $\{(\xi_n, \eta_n)\}$ where every random vector $(\xi_n, \eta_n)$ is defined on the probability space $(\Omega, \mathscr{F}, P_n)$, such that the sequence of conditional mathematical expectations $\{E(\eta_n | \xi_n = x)\}$ uniformly converges to the continuous function $f(x)$ in $[a,b]$.

This can be called random vector representation of continuous functions, which is like a bridge to be set up between real function theory and probability theory. By using this conclusion, we have a result with respect to function approximation: for any given continuous function $f(x) \in C[a,b]$, if $\{E(\eta_n | \xi_n = x)\}$ is the sequence of conditional mathematical expectations generated by the continuous function $f(x)$, then by means of $\{E(\eta_n | \xi_n = x)\}$ we can make a group of continuous base functions as follows:

$$\Phi(n) = \{\varphi_0^{(n)}(x), \varphi_1^{(n)}(x), \cdots, \varphi_n^{(n)}(x)\},$$

such that the sequence of interpolation functions formed by using $\{\Phi(n)\}$ as the following

$$f_n(x) = \sum_{l=0}^{n} \varphi_l^{(n)}(x) y_l^{(n)}, \quad n = 1, 2, 3, \cdots$$

can uniformly converge to $f(x)$.

And then, in the approximation from a sequence of random vectors to a continuous function, the base functions are appropriately selected by us, an important conclusion for quantum mechanics is deduced: classical mechanics and quantum mechanics are unified. Particularly, an interesting and very important conclusion is introduced as the fact that the mass point motion of a macroscopical object possesses a kind of wave characteristic curve, which should be called wave-mass-point duality.

## 2. The random vector presentation of continuous functions

**Lemma 2.1** Arbitrarily given $m+1$ real numbers $a_0, a_1, \cdots, a_m \in \mathbb{R}$, we denote the following symbol:

$$e_m = \max\left\{ |a_i - a_{i-1}| \,\big|\, i = 1, 2, \cdots, m \right\}.$$

And we make a permutation as the following:

$$\sigma = \begin{pmatrix} 0 & 1 & \cdots & m \\ k_0 & k_1 & \cdots & k_m \end{pmatrix},$$

such that $a_{k_0} \le a_{k_1} \le \cdots \le a_{k_m}$. If we write

$$d_m = \max\left\{ a_{k_i} - a_{k_{i-1}} \,\big|\, i = 1, 2, \cdots, m \right\},$$

Then we have that $d_m \le e_m$.

**Proof**. By the definition of $d_m$, we can know the following fact:

$$\left(\exists i \in \{1, 2, \cdots, m\}\right)\left(d_m = a_{k_i} - a_{k_{i-1}}\right).$$

If $d_m = 0$, then the conclusion of the lemma is clearly true. Now we assume that $d_m > 0$. We know that $\sigma$ is a bijection, and then $k_i \ne k_{i-1}$.

Let $s = k_i$ and $t = k_{i-1}$. So $a_s - a_t = d_m$. We consider two cases: (i) and (ii) as follows.

(i) $s < t$. If we pay attention to the total order relation:

$$a_{k_0} \le a_{k_1} \le a_{k_{i-1}} < a_{k_i} \le \cdots \le \cdots \le a_{k_m},$$

then we can learn the fact:

$$a_s, a_{s+1}, a_{s+2}, \cdots, a_{t-1}, a_t \notin \left(a_t, a_s\right) = \left(a_{k_{i-1}}, a_{k_i}\right),$$

which means $a_s, a_{s+1}, a_{s+2}, \cdots, a_{t-1}, a_t \in \left(-\infty, a_t\right] \cup \left[a_s, +\infty\right)$. Let

$$l = \min\left\{ i \in \{s, s+1, \cdots, t\} \,\big|\, a_i \in \left(-\infty, a_t\right] \right\}.$$

Clearly $l \ne s$; or else $a_l = a_s \in \left[a_s, +\infty\right)$; this will be contradictory with the fact that $a_l \in \left(-\infty, a_t\right]$. By the meaning of the subscript $l$, it is easy to understand that $a_{l-1} \in \left[a_s, +\infty\right)$. Thus we have the result:

$$e_m \ge |a_l - a_{l-1}| \ge a_s - a_t = d_m.$$

(ii) $t < s$. This time we have the following result:

$$a_t, a_{t+1}, a_{t+2}, \cdots, a_{s-1}, a_s \notin \left(a_t, a_s\right) = \left(a_{k_{i-1}}, a_{k_i}\right),$$

which means the following expression is true:

$$a_t, a_{t+1}, a_{t+2}, \cdots, a_{s-1}, a_s \in \left(-\infty, a_t\right] \cup \left[a_s, +\infty\right).$$

Let $l = \max\left\{ i \in \{t, t+1, \cdots, s\} \,\big|\, a_i \in \left(-\infty, a_t\right] \right\}$. Clearly $l \ne s$, or else we have $a_l = a_s \in \left[a_s, +\infty\right)$; this will also be contradictory with this expression $a_l \in \left(-\infty, a_t\right]$. So $a_{l+1} \in \left[a_s, +\infty\right)$. Thus we have the result:

$$e_m \ge |a_{l+1} - a_l| \ge a_s - a_t = d_m.$$

We complete the proof of the lemma. □

**Theorem 2.1** For arbitrarily given a continuous function $f(x) \in C[a, b]$, there must be a sequenc of probability spaces $\{(\Omega,\mathcal{F},P_n)\}$ and a sequence of random vectors $\{(\xi_n,\eta_n)\}$, where every random vector $(\xi_n,\eta_n)$ is defined on the probability space $(\Omega,\mathcal{F},P_n)$, such that the sequence of conditional mathematical expectations $\{E(\eta_n|\xi_n=x)\}$ uniformly converges to the continuous function $f(x)$ in $[a,b]$, that is, for any $\varepsilon>0$, there exists $N\in\mathbb{N}_+$, such that, for any $n\in\mathbb{N}_+$, if $n>N$, then

$$(\forall x\in[a,b])\left(\left|E(\eta_n|\xi_n=x)-f(x)\right|<\varepsilon\right),$$

where $\mathbb{N}_+=\{1,2,3,\cdots\}$ and $\mathbb{N}=\{0,1,2,\cdots\}$.

**Proof** **Case 1**. Let $f(x)$ be a strictly monotone functions, and we may as well assume that $f(x)$ is a strictly monotonically increasing function, as when $f(x)$ a strictly monotonically decreasing function, the proof is the same as the increasing status.

**Step 1**. Construct a group of continuous base functions.

Firstly the interval $X=[a,b]$ is equidistantly partitioned as the following:

$$a=x_0^{(n)}<x_1^{(n)}<\cdots<x_n^{(n)}=b.$$

And we write

$$X(n)=\left\{x_i^{(n)}\middle|i=0,1,\cdots,n\right\},$$

$$h(n)=\frac{b-a}{n},$$

$$x_i^{(n)}=a+ih(n),\quad i=0,1,\cdots,n.$$

Clearly $h(n)\to 0\Leftrightarrow n\to+\infty$. Then denote $Y=[c,d]=f(X)$, and put

$$y_i^{(n)}=f\left(x_i^{(n)}\right), i=0,1,\cdots,n,$$

$$Y(n)=\left\{y_i^{(n)}\middle|i=0,1,\cdots,n\right\},$$

$$c(n)=\min Y(n),$$

$$d(n)=\max Y(n)$$

Then $c(n)\ge c,\ d(n)\le d$. By using two node sets $X(n)$ and $Y(n)$, two groups of continuous base functions are formed as the following:

$$\mathcal{A}(n)=\left\{A_i^{(n)}\in C(X)\,\middle|\,i=0,1,\cdots,n\right\},$$

$$\mathcal{B}(n)=\left\{B_i^{(n)}\in C(Y)\,\middle|\,i=0,1,\cdots,n\right\}$$

where the definition of $A_i^{(n)}$ $(i=0,1,\cdots,n)$ are as the following which the figures of them are shown as Fig. 2.1 (see [1]):

$$A_0^{(n)}(x)=\begin{cases}\left(x-x_1^{(n)}\right)\Big/\left(x_0^{(n)}-x_1^{(n)}\right), & x\in\left[x_0^{(n)},x_1^{(n)}\right];\\ 0, & \text{otherwise},\end{cases}$$

$$A_i^{(n)}(x)=\begin{cases}\left(x-x_{i-1}^{(n)}\right)\Big/\left(x_i^{(n)}-x_{i-1}^{(n)}\right), & x\in\left(x_{i-1}^{(n)},x_i^{(n)}\right];\\ \left(x-x_{i+1}^{(n)}\right)\Big/\left(x_i^{(n)}-x_{i+1}^{(n)}\right), & x\in\left(x_i^{(n)},x_{i+1}^{(n)}\right];\\ 0, & \text{otherwise};\end{cases}\qquad (2.1)$$

$$i=1,2,\cdots,n-1,$$

$$A_n^{(n)}(x)=\begin{cases}\left(x-x_{n-1}^{(n)}\right)\Big/\left(x_n^{(n)}-x_{n-1}^{(n)}\right), & x\in\left(x_{n-1}^{(n)},x_n^{(n)}\right];\\ 0, & \text{otherwise},\end{cases}$$

As $X(n)=\left\{x_i^{(n)}\middle|i=0,1,\cdots,n\right\}$ is an equidistant partition node set, above $A_i^{(n)}$ $(i=0,1,\cdots,n)$ can be simplified as the following:

$$A_0^{(n)}(x)=\begin{cases}\left(x_1^{(n)}-x\right)/h(n), & x\in\left[x_0^{(n)},x_1^{(n)}\right];\\ 0\,, & \text{otherwise,}\end{cases}$$

$$A_i^{(n)}(x)=\begin{cases}\left(x-x_{i-1}^{(n)}\right)/h(n), & x\in\left(x_{i-1}^{(n)},x_i^{(n)}\right];\\ \left(x_{i+1}^{(n)}-x\right)/h(n), & x\in\left(x_i^{(n)},x_{i+1}^{(n)}\right];\\ 0, & \text{otherwise;}\end{cases}$$

$$i=1,2,\cdots,n-1,$$

$$A_n^{(n)}(x)=\begin{cases}\left(x-x_{n-1}^{(n)}\right)/h(n), & x\in\left(x_{n-1}^{(n)},x_n^{(n)}\right];\\ 0, & \text{otherwise,}\end{cases}$$

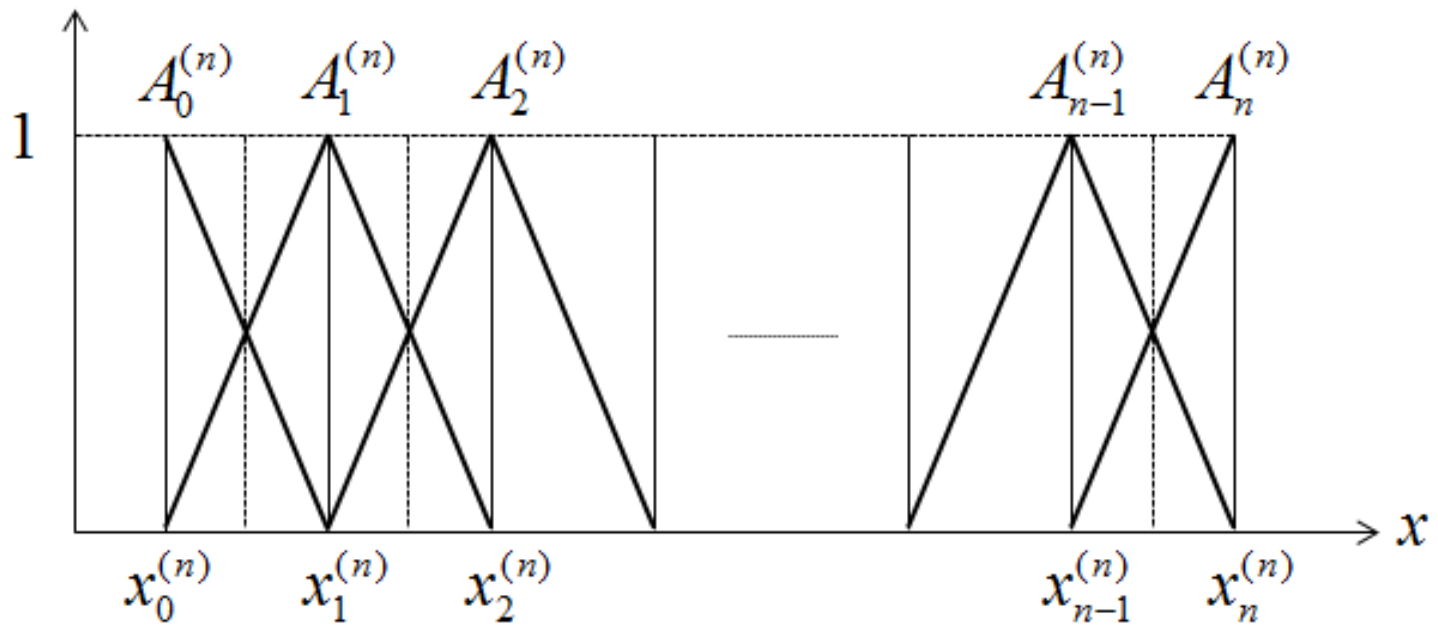

Fig. 2.1 continuous base functions $A_i^{(n)}$

As $f(x)$ being strictly monotonically increasing function, we have

$$c=c(n)=y_0^{(n)}<y_1^{(n)}<\cdots<y_n^{(n)}=d(n)=d\,. \tag{2.2}$$

By means of these nodes: $y_0^{(n)}<y_1^{(n)}<\cdots<y_n^{(n)}$, we can construct continuous base functions denoted by $B_i^{(n)}$ $(i=0,1,\cdots,n)$ which are defined on $Y$ as follows

$$B_0^{(n)}(y)=\begin{cases}\left(y-y_1^{(n)}\right)/\left(y_0^{(n)}-y_1^{(n)}\right), & y\in\left[y_0^{(n)},y_1^{(n)}\right];\\ 0\,, & \text{otherwise,}\end{cases}$$

$$B_i^{(n)}(y)=\begin{cases}\left(y-y_{i-1}^{(n)}\right)/\left(y_i^{(n)}-y_{i-1}^{(n)}\right), & y\in\left(y_{i-1},y_i\right];\\ \left(y-y_{i+1}^{(n)}\right)/\left(y_i^{(n)}-y_{i+1}^{(n)}\right), & y\in\left(y_i^{(n)}y_{i+1}^{(n)}\right];\\ 0, & \text{otherwise;}\end{cases} \tag{2.3}$$

$$i=1,2,\cdots,n-1,$$

$$B_n^{(n)}(y)=\begin{cases}\left(y-y_{n-1}^{(n)}\right)/\left(y_n^{(n)}-y_{n-1}^{(n)}\right), & y\in\left(y_{n-1}^{(n)},y_n^{(n)}\right];\\ 0, & \text{otherwise,}\end{cases}$$

So we get a group of continuous base functions $\mathcal{B}(n)=\left\{B_i^{(n)}\middle|i=0,1,\cdots,n\right\}$.

**Step 2**. Construct a sequence of probability density functions $\left\{p_n(x,y)\right\}$.

By using $\mathcal{A}(n)$ and $\mathcal{B}(n)$, we form a group of continuous base functions with two variables defined on $X\times Y=[a,b]\times[c,d]$ as follows

$$C_i^{(n)}(x,y)=A_i^{(n)}(x)\cdot B_i^{(n)}(y),$$
$$i=0,1,2,\cdots,n$$

Write $\mathcal{C}(n)=\left\{C_i^{(n)}(x,y)\middle|i=0,1,\cdots,n\right\}$, and it is easy to know that $\mathcal{C}(n)$ is a group of linearly independent functions. Then $\operatorname{span}\mathcal{C}(n)$ is just a $n+1$ dimension linear subspace of $C\left(X\times Y\right)$.

Next let $\vee=\max$, which means that for any $n+1$ real numbers $a_0,a_1,\cdots,a_n$, we have

$$\mathop{\vee}_{i=0}^{n} a_i=\max\left\{a_0,a_1,\cdots,a_n\right\}.$$

Based on $\mathcal{C}(n)$, a continuous function with two variables $R_n:X\times Y\to[0,1]$ is formed as follows

$$R_n(x,y)=\mathop{\vee}_{i=0}^{n} C_i^{(n)}(x,y)=\mathop{\vee}_{i=0}^{n}\left(A_i^{(n)}(x)\cdot B_i^{(n)}(y)\right). \tag{2.4}$$

Then we get a sequence of continuous functions with two variables $\left\{R_n(x,y)\right\}$. Then we write

$$H(n)=\int_c^d\int_a^b R_n(x,y)dxdy.$$

Clearly $\left(\forall n\in\mathbb{N}_+\right)\left(H(n)>0\right)$. And we put

$$p_n(x,y)=\frac{R_n(x,y)}{H(n)}\chi_{X\times Y}(x,y),\quad n=1,2,3,\cdots, \tag{2.5}$$

where $\chi_{X\times Y}$ is the indicative function of set $X\times Y$:

$$\chi_{X\times Y}:\mathbb{R}^2\to\{0,1\}$$
$$(x,y)\mapsto\chi_{X\times Y}(x,y)=\begin{cases}1, & (x,y)\in X\times Y,\\ 0, & (x,y)\notin X\times Y\end{cases}$$

So we get a sequence of probability density functions $\left\{p_n(x,y)\right\}$ defined on $\mathbb{R}^2=(-\infty,+\infty)^2$.

**Step 3**. Construct a sequence of probability spaces $\left\{\left(\Omega,\mathcal{F},P_n\right)\right\}$ and a sequence of random vectors $\left\{\zeta_n\right\}=\left\{\left(\xi_n,\eta_n\right)\right\}$ which every random vector $\zeta_n$ is defined on $\left(\Omega,\mathcal{F},P_n\right)$. In fact, let

$$F_n(x,y)=\int_{-\infty}^x\int_{-\infty}^y p_n(u,v)dudv;$$

then $\left\{F_n(x,y)\right\}$ is a sequence of distribution functions. Take $\Omega=\mathbb{R}^2$ and $\mathcal{F}=\mathcal{B}_2$, where $\mathcal{B}_2$ is a Borel $\sigma$ algebra on $\mathbb{R}^2$; and $P_n$ is taken as the probability measure corresponding to $F_n(x,y)$. We all know a fact in probability theory that $P_n$ must exist and be unique (see [4][9][11]). In this way we get a sequence of probability spaces as the following:

$$\left\{\left(\Omega,\mathcal{F},P_n\right)\right\}=\left\{\left(\mathbb{R}^2,\mathcal{B}_2,P_n\right)\right\}.$$

Then on every probability space $\left(\Omega,\mathcal{F},P_n\right)$ we define a random vector as follows:

$$\zeta_n=\left(\xi_n,\eta_n\right):\Omega\to\mathbb{R}^2$$
$$\omega=\left(\omega_1,\omega_2\right)\mapsto\zeta_n(\omega)=\left(\xi_n(\omega),\eta_n(\omega)\right)=\left(\omega_1,\omega_2\right)$$

For any $(x,y)\in\mathbb{R}^2$, by noticing the following expression

$$\left\{\omega\in\Omega\middle|\xi_n(\omega)\le x\right\}=\left\{\omega\in\Omega\middle|\omega_1\le x\right\}$$
$$=\left\{\omega\in\Omega\middle|\omega_1\in(-\infty,x],\omega_2\in(-\infty,+\infty)\right\}$$
$$=(-\infty,x]\times(-\infty,+\infty)\in\mathcal{B}_2=\mathcal{F}$$

We can know that $\xi_n(\omega)$ is really a random variable defined on $(\Omega, \mathcal{F}, P_n)$; in the same way, $\eta_n(\omega)$ is also really a random variable defined on $(\Omega, \mathcal{F}, P_n)$. So

$$\zeta_n(\omega) = (\xi_n(\omega), \eta_n(\omega)) = (\omega_1, \omega_2) = \omega$$

is just a random vector defined on $(\Omega, \mathcal{F}, P_n)$.

Let $F_{\zeta_n}(x, y)$ be the distribution function of the random vector $\zeta_n(\omega)$. For any $(x, y) \in \mathbb{R}^2$, because $P_n$ is the probability measure corresponding to $F_n(x, y)$, we have

$$\begin{aligned} F_{\zeta_n}(x, y) &= P_n\left(\left\{\omega \in \Omega \middle| \xi_n(\omega) \le x, \eta_n(\omega) \le y\right\}\right) \\ &= P_n\left(\left\{\omega \in \Omega \middle| \omega_1 \le x, \omega_2 \le y\right\}\right) \\ &= P_n\left((-\infty, x] \times (-\infty, y]\right) = F_n(x, y) \end{aligned}$$

i.e., the distribution function $F_{\zeta_n}(x, y)$ of the random vector $\zeta_n(\omega)$ is just $F_n(x, y)$, which means

$$F_{\zeta_n}(x, y) \equiv F_n(x, y).$$

So the sequence of conditional expectations of the sequence of random vectors $\{\zeta_n\} = \{(\xi_n, \eta_n)\}$ is just $\{E(\eta_n | \xi_n = x)\}$, where

$$\begin{aligned} E(\eta_n | \xi_n = x) &= \frac{\int_{-\infty}^{+\infty} y p_n(x, y) dy}{\int_{-\infty}^{+\infty} p_n(x, y) dy} \\ &= \frac{\int_c^d y R_n(x, y) dy}{\int_c^d R_n(x, y) dy}, \quad n = 1, 2, 3, \cdots \end{aligned}$$

**Step 4**. Prove the fact that the sequence of conditional expectations $\{E(\eta_n | \xi_n = x)\}$ converges to $f(x)$ everywhere in $[a, b]$. Firstly we prove a result as the following:

$$(\forall n \in \mathbb{N}_+)(\forall x \in [a, b])\left(\int_c^d R_n(x, y) dy > 0\right).$$

In fact, for any $x \in [a, b]$, clearly $(\exists i \in \{1, 2, \cdots, n\})\left(x \in \left[x_{i-1}^{(n)}, x_i^{(n)}\right]\right)$. Then for any $y \in [c, d]$, we have

$$\begin{aligned} R_n(x, y) &= \bigvee_{k=0}^{n} \left(A_k^{(n)}(x) \cdot B_k^{(n)}(y)\right) \\ &= \begin{cases} \left(A_{i-1}^{(n)}(x) \cdot B_{i-1}^{(n)}(y)\right) \vee \left(A_i^{(n)}(x) \cdot B_i^{(n)}(y)\right), & y \in \left[y_{i-2}^{(n)}, y_{i+1}^{(n)}\right], \\ 0, & y \in [c, d] - \left[y_{i-2}^{(n)}, y_{i+1}^{(n)}\right] \end{cases} \end{aligned}$$

It is easy to learn the fact that

$$\left(\exists y' \in \left(y_{i-2}^{(n)}, y_{i+1}^{(n)}\right)\right)\left(R_n(x, y') > 0\right).$$

For this fixed $x \in [a, b]$, since $R_n(x, y)$ is a continuous function with respect to $y$, we have

$$(\exists \delta > 0)\left((y' - \delta, y' + \delta) \subset \left(y_{i-2}^{(n)}, y_{i+1}^{(n)}\right)\right)$$

$$\left(\forall y \in (y' - \delta, y' + \delta)\right)\left(R_n(x, y) > 0\right)$$

By means of mean value theorem of integrals, there exists a point $\xi \in \left[y' - \frac{\delta}{2}, y' + \frac{\delta}{2}\right]$, such that

$$\int_c^d R_n(x, y) dy \ge \int_{(y' - \delta, y' + \delta)} R_n(x, y) dy$$

$$\geq \int_{\left[y'-\frac{\delta}{2},y'+\frac{\delta}{2}\right]} R_n(x,y)dy = R_n(x,\xi)\cdot\delta > 0$$

So above result is correct. By this result, for any $n \in \mathbb{N}_+$, the following expression

$$E\left(\eta_n \middle| \xi_n = x\right) = \frac{\int_c^d yR_n(x,y)dy}{\int_c^d R_n(x,y)dy}$$

is meaningful.

Now we turn to prove the fact that the sequence of unary functions $\left\{E\left(\eta_n \middle| \xi_n = x\right)\right\}$ converges to $f(x)$ everywhere in $[a,b]$. In fact, for any $x \in [a,b]$, clearly we have

$$\left(\exists i \in \{1,2,\cdots,n\}\right)\left(x \in \left[x_{i-1}^{(n)}, x_i^{(n)}\right]\right).$$

Then we can get the following expression:

$$\begin{aligned} R_n(x,y) &= \bigvee_{k=0}^{n}\left(A_k^{(n)}(x)\cdot B_k^{(n)}(y)\right) \\ &= \begin{cases} \left(A_{i-1}^{(n)}(x)\cdot B_{i-1}^{(n)}(y)\right)\vee\left(A_i^{(n)}(x)\cdot B_i^{(n)}(y)\right), & y \in \left[y_{i-2}^{(n)}, y_{i+1}^{(n)}\right], \\ 0, & y \in [c,d]-\left[y_{i-2}^{(n)}, y_{i+1}^{(n)}\right] \end{cases} \end{aligned}$$

By means of the first mean value theorem of integrals, there exsits a point $\eta_n(x) \in \left[y_{i-2}^{(n)}, y_{i+1}^{(n)}\right]$, such that

$$\begin{aligned} E\left(\eta_n \middle| \xi_n = x\right) &= \frac{\int_c^d yR_n(x,y)dy}{\int_c^d R_n(x,y)dy} \\ &= \frac{\int_{y_{i-2}^{(n)}}^{y_{i+1}^{(n)}} yR_n(x,y)dy}{\int_{y_{i-2}^{(n)}}^{y_{i+1}^{(n)}} R_n(x,y)dy} = \frac{\eta_n(x)\int_{y_{i-2}^{(n)}}^{y_{i+1}^{(n)}} R_n(x,y)dy}{\int_{y_{i-2}^{(n)}}^{y_{i+1}^{(n)}} R_n(x,y)dy} \\ &= \eta_n(x) \end{aligned}$$

For $f(x)$ being continuous and by noticing that $\eta_n(x) \in \left[y_{i-2}^{(n)}, y_{i+1}^{(n)}\right]$, based on the medium value theorem for continuous functions, we have

$$\left(\exists \overline{x} \in \left[x_{i-2}^{(n)}, x_{i+1}^{(n)}\right]\right)\left(f\left(\overline{x}\right) = \eta_n(x)\right).$$

Since $f(x) \in \left[y_{i-1}^{(n)}, y_{i+1}^{(n)}\right] \subset \left[y_{i-2}^{(n)}, y_{i+1}^{(n)}\right]$, we know that

$$n \to \infty \Rightarrow \left|y_{i+1}^{(n)} - y_{i-2}^{(n)}\right| \to 0 \Rightarrow \eta_n(x) = f\left(\overline{x}\right) \to f(x). \tag{2.6}$$

For $x \in [a,b]$ being arbitrary, we get the fact as follows

$$\left(\forall x \in [a,b]\right)\left(\lim_{n\to\infty} E\left(\eta_n \middle| \xi_n = x\right) = f(x)\right).$$

**Step 5**. Prove the fact that the sequence of conditional expectations $\left\{E\left(\eta_n \middle| \xi_n = x\right)\right\}$ uniformly converges to $f(x)$ in $[a,b]$.

Because $f(x)$ is continuous in $[a,b]$, for any $\varepsilon > 0$, there exists $N \in \mathbb{N}_+$, such that

$$\left(\forall n \in \mathbb{N}_+\right)\left(n > N \Rightarrow \max\left\{\left|\Delta y_i^{(n)}\right| \middle| i = 1,2,\cdots,n\right\} < \frac{\varepsilon}{3}\right)$$

By (2.6), for any $x \in [a,b]$, when $n > N$, we have

$$\left|E\left(\eta_n \middle| \xi_n = x\right) - f(x)\right| = \left|\eta_n(x) - f(x)\right|$$

$$\le \left| y_{i+1}^{(n)} - y_{i-2}^{(n)} \right| \le 3\max\left\{ \left| \Delta y_i^{(n)} \right| \middle| i = 1,2,\cdots,n \right\}$$

$$< 3\cdot\frac{\varepsilon}{3} = \varepsilon$$

So far we have proved the conclusion: for any $\varepsilon > 0$, there exists $N \in \mathbb{N}_+$, such that, for any $n \in \mathbb{N}_+$, if $n > N$, then

$$(\forall x \in [a,b])\left(\left| E\left(\eta_n \middle| \xi_n = x\right) - f(x) \right| < \varepsilon\right)$$

This means that the sequence of conditional expectations $\left\{ E\left(\eta_n \middle| \xi_n = x\right) \right\}$ uniformly converges to $f(x)$ in $[a,b]$.

**Case 2**. $f(x)$ is not a strictly monotonic function and not a constant function.

This time, the elements in $Y(n)$ are not of monotonic property like as $y_0 \le y_1 \le \cdots \le y_n$ with respect to the subscripts $i$, which brings some difficulty to construct continuous base functions $B_i$ being shaped like (2.3); so we should make a permutation on the subscript set $\{0,1,\cdots,n\}$ as the following:

$$\sigma = \begin{pmatrix} 0 & 1 & \cdots & n \\ k_0 & k_1 & \cdots & k_n \end{pmatrix},$$

$$\left(\forall i \in \{0,1,\cdots,n\}\right)\left(k_i = \sigma(i)\right)$$

such that the permutated subscript set $K(n) = \left\{k_0, k_1, \cdots, k_n\right\}$ is with the condition:

$$c(n) = y_{k_0}^{(n)} \le y_{k_1}^{(n)} \le \cdots \le y_{k_n}^{(n)} = d(n). \tag{2.7}$$

As (2.7) is not with strictly monotonicity about the subscript, we will deal with it in two situations.

1) Assume $c(n) = y_{k_0}^{(n)} < y_{k_1}^{(n)} < \cdots < y_{k_n}^{(n)} = d(n)$. By using these nodes

$$y_{k_0}^{(n)}, y_{k_1}^{(n)}, \cdots, y_{k_n}^{(n)}$$

in $[c,d]$, we form continuous base functions $B_{k_j}^{(n)}$ $\left(j = 0,1,\cdots,n\right)$ as follows

$$B_{k_0}^{(n)}(y) = \begin{cases} \left(y - y_{k_1}^{(n)}\right) \big/ \left(y_{k_0}^{(n)} - y_{k_1}^{(n)}\right), & y \in \left[y_{k_0}^{(n)}, y_{k_1}^{(n)}\right]; \\ 0, & \text{otherwise}, \end{cases}$$

$$B_{k_j}^{(n)}(y) = \begin{cases} \left(y - y_{k_{j-1}}^{(n)}\right) \big/ \left(y_{k_j}^{(n)} - y_{k_{j-1}}^{(n)}\right), & y \in \left(y_{k_{j-1}}^{(n)}, y_{k_j}^{(n)}\right]; \\ \left(y - y_{k_{j+1}}^{(n)}\right) \big/ \left(y_{k_j}^{(n)} - y_{k_{j+1}}^{(n)}\right), & y \in \left(y_{k_j}^{(n)}, y_{k_{j+1}}^{(n)}\right]; \\ 0, & \text{otherwise}; \end{cases}$$

$$j = 1,2,\cdots,n-1,$$

$$B_{k_n}^{(n)}(y) = \begin{cases} \left(y - y_{k_{n-1}}^{(n)}\right) \big/ \left(y_{k_n}^{(n)} - y_{k_{n-1}}^{(n)}\right), & y \in \left(y_{k_{n-1}}^{(n)}, y_{k_n}^{(n)}\right]; \\ 0, & \text{otherwise}, \end{cases}$$

And we get a sequence of continuous functions with two variables $\left\{R_n(x,y)\right\}$ as the following

$$R_n(x,y) = \bigvee_{j=0}^{n} C_{k_j}^{(n)}(x,y)$$

$$= \bigvee_{j=0}^{n} \left(A_{k_j}^{(n)}(x) \cdot B_{k_j}^{(n)}(y)\right)$$

Like Case 1, we can get a sequence of probability density functions $\left\{p_n(x,y)\right\}$ defined on the real number field $\mathbb{R}^2 = (-\infty, +\infty)^2$ and a sequence of random vectors $\left\{\zeta_n\right\}$ $= \left\{\left(\xi_n, \eta_n\right)\right\}$; then we have a sequence of conditional mathematical expectations $\left\{E\left(\eta_n \middle| \xi_n = x\right)\right\}$, where

$$E\left(\eta_n \middle| \xi_n = x\right) = \frac{\int_{-\infty}^{+\infty} y p_n(x,y)dy}{\int_{-\infty}^{+\infty} p_n(x,y)dy}$$

$$= \frac{\int_c^d yR_n(x,y)dy}{\int_c^d R_n(x,y)dy}, \quad n = 1,2,3,\cdots$$

Now we prove the fact that the sequence of conditional mathematical expectations $\left\{E\left(\eta_n \middle| \xi_n = x\right)\right\}$ uniformly converges to $f(x)$ in $[a,b]$.

Firstly it is easy to know that the fact as the following:

$$\left(\forall n \in \mathbb{N}_+\right)\left(\forall x \in [a,b]\right)\left(\int_c^d R_n(x,y)dy > 0\right).$$

Here we should stipulate that $k_{-1} = k_0$ and $k_{n+1} = k_n$. For any $x \in [a,b]$, clearly we have

$$\left(\exists s,t \in \{0,1,\cdots,n\}\right)\left(x \in \left[x_{k_s}, x_{k_t}\right]\right).$$

Then we get the following expression:

$$R_n(x,y) = \bigvee_{j=0}^{n}\left(A_{k_j}^{(n)}(x)\cdot B_{k_j}^{(n)}(y)\right)$$

$$= \begin{cases} \left(A_{k_s}^{(n)}(x)\cdot B_{k_s}^{(n)}(y)\right) \vee \left(A_{k_t}^{(n)}(x)\cdot B_{k_t}^{(n)}(y)\right), & y \in \left[y_{k_{s-1}}^{(n)}, y_{k_{s+1}}^{(n)}\right] \cup \left[y_{k_{t-1}}^{(n)}, y_{k_{t+1}}^{(n)}\right], \\ 0, & y \in [c,d] - \left(\left[y_{k_{s-1}}^{(n)}, y_{k_{s+1}}^{(n)}\right] \cup \left[y_{k_{t-1}}^{(n)}, y_{k_{t+1}}^{(n)}\right]\right) \end{cases}$$

Let $y_* = \min\left\{y_{k_{s-1}}^{(n)}, y_{k_{t-1}}^{(n)}, y_{k_{s+1}}^{(n)}, y_{k_{t+1}}^{(n)}\right\}$ and $y^* = \max\left\{y_{k_{s-1}}^{(n)}, y_{k_{t-1}}^{(n)}, y_{k_{s+1}}^{(n)}, y_{k_{t+1}}^{(n)}\right\}$. It assumes that $y_* = y_{k_{s-1}}^{(n)}$ and $y^* = y_{k_{t+1}}^{(n)}$, and clearly we know the fact:

$$\left[y_*, y^*\right] = \left[y_{k_{s-1}}^{(n)}, y_{k_{t+1}}^{(n)}\right] \supset \left[y_{k_{s-1}}^{(n)}, y_{k_{s+1}}^{(n)}\right] \cup \left[y_{k_{t-1}}^{(n)}, y_{k_{t+1}}^{(n)}\right].$$

So above expression can be written as the following:

$$R_n(x,y) = \bigvee_{j=0}^{n}\left(A_{k_j}^{(n)}(x)\cdot B_{k_j}^{(n)}(y)\right)$$

$$= \begin{cases} \left(A_{k_s}^{(n)}(x)\cdot B_{k_s}^{(n)}(y)\right) \vee \left(A_{k_t}^{(n)}(x)\cdot B_{k_t}^{(n)}(y)\right), & y \in \left[y_{k_{s-1}}^{(n)}, y_{k_{t+1}}^{(n)}\right], \\ 0, & y \in [c,d] - \left[y_{k_{s-1}}^{(n)}, y_{k_{t+1}}^{(n)}\right] \end{cases}$$

By first means of mean value theorem for integrals, there exsits a point $\eta_n(x) \in \left[y_{k_{s-1}}^{(n)}, y_{k_{t+1}}^{(n)}\right]$, such that

$$E\left(\eta_n \middle| \xi_n = x\right) = \frac{\int_c^d yR_n(x,y)dy}{\int_c^d R_n(x,y)dy}$$

$$= \frac{\int_{y_{k_{s-1}}^{(n)}}^{y_{k_{t+1}}^{(n)}} yR_n(x,y)dy}{\int_{y_{k_{s-1}}^{(n)}}^{y_{k_{t+1}}^{(n)}} R_n(x,y)dy} = \frac{\eta_n(x)\int_{y_{k_{s-1}}^{(n)}}^{y_{k_{t+1}}^{(n)}} R_n(x,y)dy}{\int_{y_{k_{s-1}}^{(n)}}^{y_{k_{t+1}}^{(n)}} R_n(x,y)dy}$$

$$= \eta_n(x)$$

Write $d_n = \max\left\{\Delta y_{k_j}^{(n)} \middle| j = 1,2,\cdots,n\right\}$, where $\Delta y_{k_j}^{(n)} = y_{k_j}^{(n)} - y_{k_{j-1}}^{(n)}$, $j = 1,2,\cdots,n$. We can prove the fact that $\lim\limits_{n\to\infty} d_n = 0$. In fact, let

$$e_n = \max\left\{\left|y_i^{(n)} - y_{i-1}^{(n)}\right| \middle| i = 1, 2, \cdots, n\right\},$$

by means of lemma 2.1 we have $(\forall n \in \mathbb{N}_+)(d_n \le e_n)$. By using this result, because $f(x)$ is uniformly continuous in $[a,b]$, it is true that $\lim_{n\to\infty} e_n = 0$; so we have $\lim_{n\to\infty} d_n = 0$.

And then by above result we have the following expression:

$$(\forall \varepsilon > 0)(\exists N_1 \in \mathbb{N}_+)(\forall n \in \mathbb{N}_+)\left(n > N_1 \Rightarrow d_n < \frac{\varepsilon}{3}\right).$$

By noticing the fact that

$$n \to \infty \Rightarrow x_{k_s}^{(n)} - x_{k_t}^{(n)} \to 0$$

$$\Rightarrow \left|y_{k_s}^{(n)} - y_{k_t}^{(n)}\right| = \left|f\left(x_{k_s}^{(n)}\right) - f\left(x_{k_t}^{(n)}\right)\right| \to 0$$

We immediately have the result: there exists $N_2 \in \mathbb{N}_+$, for any $n \in \mathbb{N}_+$, such that

$$n > N_2 \Rightarrow \left|y_{k_s}^{(n)} - y_{k_t}^{(n)}\right| < \frac{\varepsilon}{3}. \tag{2.8}$$

Take $N = \max\{N_1, N_2\}$, and when $n > N$, we get the following expression:

$$\left|y_{k_{t+1}}^{(n)} - y_{k_{s-1}}^{(n)}\right| \le \left|y_{k_{t+1}}^{(n)} - y_{k_t}^{(n)}\right| + \left|y_{k_t}^{(n)} - y_{k_s}^{(n)}\right| + \left|y_{k_s}^{(n)} - y_{k_{s-1}}^{(n)}\right|$$

$$< d_n + \frac{\varepsilon}{3} + d_n < \frac{\varepsilon}{3} + \frac{\varepsilon}{3} + \frac{\varepsilon}{3} = \varepsilon$$

By this result we know that $\left|y_{k_{t+1}}^{(n)} - y_{k_{s-1}}^{(n)}\right| \xrightarrow{n\to\infty} 0$, so we get

$$\left|y_{k_s}^{(n)} - \eta_n(x)\right| \xrightarrow{n\to\infty} 0.$$

By means of (2.8), we can have $\left|f(x) - y_{k_s}^{(n)}\right| \xrightarrow{n\to\infty} 0$. At last we have

$$\left|f(x) - \eta_n(x)\right| \le \left|f(x) - y_{k_s}^{(n)}\right| + \left|y_{k_s}^{(n)} - \eta_n(x)\right| \xrightarrow{n\to\infty} 0.$$

i.e., $\eta_n(x) \xrightarrow{n\to\infty} f(x)$. So we get the conclusion:

$$(\forall x \in [a,b])\left(\lim_{n\to\infty} E\left(\eta_n \middle| \xi_n = x\right) = f(x)\right).$$

And then, similar to above proof, we can obtain the conclusion that the sequence of conditional mathematical expectations $\left\{E\left(\eta_n \middle| \xi_n = x\right)\right\}$ uniformly converges to $f(x)$ in $[a,b]$.

2) Assume $c(n) = y_{k_0}^{(n)} \le y_{k_1}^{(n)} \le \cdots \le y_{k_n}^{(n)} = d(n)$. The elements in the set

$$Y(n) = \left\{y_{k_0}^{(n)}, y_{k_1}^{(n)}, \cdots, y_{k_n}^{(n)}\right\}$$

should be screened firstly. In fact, write $K(n) = \{k_0, k_1, \cdots, k_n\}$, and an equivalent relation $\sim$ is defined as the following:

$$(\forall s,t \in \{0,1,\cdots,n\})\left(k_s \sim k_t \Leftrightarrow y_{k_s} = y_{k_t}\right).$$

We can get a quotient set of $K(n)$ as being

$$K(n)/\sim = \left\{[k_j] \middle| j = 0,1,\cdots,n\right\},$$

where $[k_j]$ is the equivalence class what $k_j$ belongs to. Suppose the elements of $K(n)/\sim$ are as the following:

$$[k_{j_0}], [k_{j_1}], \cdots, [k_{j_{q(n)}}],$$

where $0 \le q(n) \le n$, and we stipulate that the representation element $k_{j_s}$ is taken as the least element in $\left[k_{j_s}\right]$. Then we have

$$y_{k_{j0}}^{(n)} < y_{k_{j1}}^{(n)} < \cdots < y_{k_{jq(n)}}^{(n)} .$$

By using the nodes $y_{k_{j0}}^{(n)}, y_{k_{j1}}^{(n)}, \cdots, y_{k_{jq(n)}}^{(n)}$ in $[c,d]$, a group of continuous base functions

$$B_{k_{j_s}}^{(n)}, \quad s = 0,1,\cdots,q(n)$$

are formed as the following:

$$B_{k_{j0}}^{(n)}(y) = \begin{cases} \left(y - y_{k_{j1}}^{(n)}\right) \Big/ \left(y_{k_{j0}}^{(n)} - y_{k_{j1}}^{(n)}\right), & y \in \left[y_{k_{j0}}^{(n)}, y_{k_{j1}}^{(n)}\right]; \\ 0\,, & \text{otherwise}, \end{cases}$$

$$B_{k_{j_s}}^{(n)}(y) = \begin{cases} \left(y - y_{k_{j_{s-1}}}^{(n)}\right) \Big/ \left(y_{k_{j_s}}^{(n)} - y_{k_{j_{s-1}}}^{(n)}\right), & y \in \left(y_{k_{j_{s-1}}}^{(n)}, y_{k_{j_s}}^{(n)}\right]; \\ \left(y - y_{k_{j_{s+1}}}^{(n)}\right) \Big/ \left(y_{k_{j_s}}^{(n)} - y_{k_{j_{s+1}}}^{(n)}\right), & y \in \left(y_{k_{j_s}}^{(n)}, y_{k_{j_{s+1}}}^{(n)}\right]; \\ 0, & \text{otherwise}; \end{cases} \tag{2.9}$$

$$s = 1,2,\cdots,q(n)-1,$$

$$B_{k_{jq(n)}}^{(n)}(y) = \begin{cases} \left(y - y_{k_{jq(n)-1}}^{(n)}\right) \Big/ \left(y_{k_{jq(n)}}^{(n)} - y_{k_{jq(n)-1}}^{(n)}\right), & y \in \left(y_{k_{jq(n)-1}}, y_{k_{jq(n)}}^{(n)}\right]; \\ 0, & \text{otherwise}, \end{cases}$$

So at the nodes $y_{k_{j0}}^{(n)}, y_{k_{j1}}^{(n)}, \cdots, y_{k_{jq(n)}}^{(n)}$ with respect to the representation elements $k_{j_0}, k_{j_1}, \cdots, k_{j_{q(n)}}$, the group of continuous base functions asfollows

$$B_{k_{j0}}^{(n)}(y), B_{k_{j1}}^{(n)}(y), \cdots, B_{k_{jq(n)}}^{(n)}(y)$$

has been defined.

And then, we should stipulate: for every equivalence class $\left[k_{j_s}\right]$, all the elements in $\left[k_{j_s}\right]$ have been-corresponded to the same continuous base function $B_{k_{j_s}}^{(n)}(y)$. Therefore, all the nodes in $[c,d]$,

$$y_{k_0}^{(n)} \le y_{k_1}^{(n)} \le \cdots \le y_{k_n}^{(n)},$$

have been defined their continuous base functions: $B_{k_0}^{(n)}(y), B_{k_1}^{(n)}(y), \cdots, B_{k_n}^{(n)}(y)$. Hence, we can get a sequence of continuous base functions with two variables $\left\{R_n(x,y)\right\}$ as follows

$$\begin{aligned} R_n(x,y) &= \mathop{\vee}_{j=0}^{n} C_{k_j}^{(n)}(x,y) \\ &= \mathop{\vee}_{j=0}^{n} \left(A_{k_j}^{(n)}(x) \cdot B_{k_j}^{(n)}(y)\right), \\ & n = 1,2,3,\cdots \end{aligned} \tag{2.10}$$

Similar to above method we have ever been used, we can have a sequence of probability density functions $\left\{p_n(x,y)\right\}$ defined on $\mathbb{R}^2 = (-\infty,+\infty)^2$ and a sequence of random vectors $\left\{\zeta_n\right\} = \left\{\left(\xi_n, \eta_n\right)\right\}$; then we immediately get a sequence of conditional mathematical expectations $\left\{E\left(\eta_n \middle| \xi_n = x\right)\right\}$, where

$$\begin{aligned} E\left(\eta_n \middle| \xi_n = x\right) &= \frac{\int_{-\infty}^{+\infty} y p_n(x,y)dy}{\int_{-\infty}^{+\infty} p_n(x,y)dy} \\ &= \frac{\int_c^d y R_n(x,y)dy}{\int_c^d R_n(x,y)dy}, \quad n = 1,2,3,\cdots \end{aligned}$$

Now we should prove the fact that the sequence of conditional mathematical expectations $\left\{E\left(\eta_n \middle| \xi_n = x\right)\right\}$ uniformly converges to $f(x)$ in $[a,b]$.

Firstly it is not difficult to know the fact that

$$\left(\forall n \in \mathbb{N}_+\right)\left(\forall x \in [a,b]\right)\left(\int_c^d R_n(x,y)dy > 0\right).$$

For any fixed $x \in [a,b]$, $\left(\exists s,t \in \{0,1,\cdots,n\}\right)\left(x \in \left[x_{k_s}^{(n)}, x_{k_t}^{(n)}\right]\right)$.

Next we should consider the following two situations.

i) When $B_{k_s}^{(n)}(y) \equiv B_{k_t}^{(n)}(y)$, we have

$$\begin{aligned} R_n(x,y) &= \bigvee_{j=0}^{n}\left(A_{k_j}^{(n)}(x)\cdot B_{k_j}^{(n)}(y)\right) \\ &= \begin{cases} \left(A_{k_s}^{(n)}(x)\cdot B_{k_s}^{(n)}(y)\right)\vee\left(A_{k_t}^{(n)}(x)\cdot B_{k_t}^{(n)}(y)\right), & y \in \left[y_{k_{s-1}}^{(n)}, y_{k_{s+1}}^{(n)}\right], \\ 0, & y \in [c,d] - \left[y_{k_{s-1}}^{(n)}, y_{k_{s+1}}^{(n)}\right] \end{cases} \end{aligned}$$

ii) When 当 $B_{k_s}^{(n)}(y) \not\equiv B_{k_t}^{(n)}(y)$, we should have

$$\begin{aligned} R_n(x,y) &= \bigvee_{j=0}^{n}\left(A_{k_j}^{(n)}(x)\cdot B_{k_j}^{(n)}(y)\right) \\ &= \begin{cases} \left(A_{k_s}^{(n)}(x)\cdot B_{k_s}^{(n)}(y)\right)\vee\left(A_{k_t}^{(n)}(x)\cdot B_{k_t}^{(n)}(y)\right), & y \in \left[y_{k_{s-1}}^{(n)}, y_{k_{s+1}}^{(n)}\right]\cup\left[y_{k_{t-1}}^{(n)}, y_{k_{t+1}}^{(n)}\right], \\ 0, & y \in [c,d] - \left(\left[y_{k_{s-1}}^{(n)}, y_{k_{s+1}}^{(n)}\right]\cup\left[y_{k_{t-1}}^{(n)}, y_{k_{t+1}}^{(n)}\right]\right) \end{cases} \end{aligned}$$

However, in either case, we can use the method similar to 1) to prove the result:

$$\left(\forall x \in [a,b]\right)\left(\lim_{n\to\infty} E\left(\eta_n \middle| \xi_n = x\right) = f(x)\right),$$

and the sequence of conditional mathematical expectations $\left\{E\left(\eta_n \middle| \xi_n = x\right)\right\}$ uniformly converges to $f(x)$ in $[a,b]$.

**Case 3**. $f(x)$ is a constant function, i.e., $\left(\exists \beta \in \mathbb{R}\right)\left(\forall x \in [a,b]\right)\left(f(x) = \beta\right)$. Clearly this is a kind of degrading situation. So we should take a distribution function:

$$F_\eta(y) = \begin{cases} 0, & y \in (-\infty, \beta), \\ 1, & y \in [\beta, +\infty) \end{cases}$$

And we construct a probability space $\left(\Omega, \mathcal{F}, P\right)$, where $P$ is a probability measure corresponding to $F_\eta(y)$, $\Omega = \mathbb{R}^1$ and $\mathcal{F} = \mathcal{B}_1$. Take the random variable as follows

$$\eta : \Omega \to \mathbb{R}^1, \omega \mapsto \eta(\omega) = \omega.$$

It is easy to know that the distribution function of $\eta$ is just $F_\eta(y)$. By noticing the following expression

$$\begin{aligned} &P\left(\left\{\omega \in \Omega \middle| \eta(\omega) = \beta\right\}\right) = P\left(\left\{\omega \in \Omega \middle| \omega = \beta\right\}\right) \\ &= F_\eta(\beta) - F_\eta(\beta - 0) = 1 - 0 = 1 \end{aligned}$$

We know that $E(\eta) = \beta$. Hence $f(x) \equiv E(\eta)$.

Of course, we can also take another random variable $\xi$ defined on $\left(\Omega, \mathcal{F}, P\right)$, which $\xi$ is required to be independent with $\eta$. So $(\xi, \eta)$ can be regarded as a random vector on $\left(\Omega, \mathcal{F}, P\right)$. And we have

$$E\left(\eta \middle| \xi = x\right) \equiv E(\eta) = \beta.$$

Furthermore, we tale a sequence of random vectors $\left\{\left(\xi_n, \eta_n\right)\right\}$, such that

$$\left(\forall n \in \mathbb{N}_+\right)\left(\left(\xi_n, \eta_n\right) = (\xi, \eta)\right),$$

Then $\left\{E\left(\eta_n \middle| \xi_n = x\right)\right\}$ uniformly converges to $f(x)$ in $[a,b]$. We finish the proof of theorem 2.1. □

## 3. The Significance of Function Approximation of Theorem 2.1

In above section, we have proved the conclusion: the sequence of conditional mathematical expectations $\left\{E\left(\eta_n \middle| \xi_n = x\right)\right\}$ uniformly converges to $f(x) \in C[a,b]$ in $[a,b]$. Now we reveal the significance of function approximation of $\left\{E\left(\eta_n \middle| \xi_n = x\right)\right\}$ about continuous function $f(x)$. We only consider the continuous function space $C[a,b]$. In $C[a,b]$, addition operation "$+$" and scalar multiplication operation "$\cdot$" are defined as the following:

$$\begin{aligned}
&+ : C[a,b]\times C[a,b] \to C[a,b] \\
&\qquad (f,g) \mapsto +(f,g) = f+g, \\
&\left(\forall x\in[a,b]\right)\left[(f+g)(x) = f(x)+g(x)\right]; \\
&\cdot : \mathbb{R}\times C[a,b] \to C[a,b] \\
&\qquad (a,f) \mapsto \cdot(a,f) = a\cdot f, \\
&\left(\forall x\in[a,b]\right)\left[(a\cdot f)(x) = a\cdot f(x)\right]
\end{aligned}$$

We all know that $\left(C[a,b],+,\mathbb{R},\cdot\right)$ forms a linear space, which can be simply denoted by $C[a,b]$. In $C[a,b]$, we define a norm operation as follows

$$\begin{aligned}
&\|\cdot\| : C[a,b] \to [0,+\infty) \\
&f \mapsto \|\cdot\|(f) = \|f\| = \max_{x\in[a,b]}|f(x)|
\end{aligned}$$

Then $\left(C[a,b],\|\cdot\|\right)$ is a normed linear space, which can be also denoted by $C[a,b]$ (see [6][7][10]).

Clearly $C[a,b]$ is an infinite dimension normed linear space. Suppose $f(x)\in C[a,b]$ is a "complicated" function; for every $n\in\mathbb{N}_+$, we try to find a group of $n+1$ linearly independent "simple" functions as following:

$$\Phi(n) = \left\{\varphi_0^{(n)}(x), \varphi_1^{(n)}(x), \cdots, \varphi_n^{(n)}(x)\right\} \subset C[a,b],$$

and $n+1$ real numbers $a_0^{(n)}, a_1^{(n)}, \cdots, a_n^{(n)} \in \mathbb{R}$, where there at least exists one real number $a_i^{(n)} \neq 0$, such that $f(x)$ can be approximately expressed by $\Phi(n)$, i.e.,

$$\left(\forall x\in[a,b]\right)\left(\left|f(x)-\sum_{i=0}^{n} a_i^{(n)}\varphi_i^{(n)}(x)\right| < \varepsilon\right). \tag{3.1}$$

where $\varepsilon > 0$ is a kind of approximation accuracy determined in advance. If we put

$$f_n(x) = \sum_{i=0}^{n} a_i^{(n)}\varphi_i^{(n)}(x),$$

then we have a sequence of continuous functions $\left\{f_n(x)\right\}_{n=1}^{\infty}$ in $C[a,b]$. Expression (3.1) means that the sequence of continuous functions $\left\{f_n(x)\right\}_{n=1}^{\infty}$ uniformly converges to $f(x)$ in $[a,b]$, i.e., for any $\varepsilon>0$, there exists $N\in\mathbb{N}_+$, such that

$$\left(\forall n\in\mathbb{N}_+\right)\left(n>N \Rightarrow \|f_n - f\| < \varepsilon\right).$$

For some fixed $\varepsilon$, when $n>N$, $\left(\mathrm{span}\Phi(n),\|\cdot\|\right)$ is a $n+1$ dimension normed linear sunspace of $\left(C[a,b],\|\cdot\|\right)$, such that, there exist real numbers $a_0^{(n)}, a_1^{(n)}, \cdots, a_n^{(n)} \in \mathbb{R}$, such that $\|f_n - f\| < \varepsilon$, where $\mathrm{span}\Phi(n)$ means that a normed linear subspace of $C[a,b]$ generated by $\Phi(n)$. In other words, on this $\varepsilon$, we can use a kind of linear combination of the base functions in $\mathrm{span}\Phi(n)$ as follows

$$f_n(x) = \sum_{i=1}^{n} a_i^{(n)} \varphi_i^{(n)}(x)$$

to take the place of $f(x)$ approximately, or we say that $f_n(x)$ can approximate $f(x)$ to the approximation accuracy $\varepsilon$. This is one of the basic ideas of function approximation.

Particularly, function interpolation is a kind of function approximation method commonly used by us (see [12][13]). Based on this idea, firstly the interval $[a,b]$ is partitioned as

$$a = x_0^{(n)} < x_1^{(n)} < \cdots < x_n^{(n)} = b,$$

where the partition may not be equidistant. Write

$$X(n) = \left\{ x_i^{(n)} \middle| i = 0,1,\cdots,n \right\},$$

$$y_i^{(n)} = f\left(x_i^{(n)}\right), i = 0,1,\cdots,n,$$

$$Y(n) = \left\{ y_i^{(n)} \middle| i = 0,1,\cdots,n \right\}$$

By the use of the node set $X(n)$, a group of base function as the following

$$\Phi(n) = \left\{ \varphi_0^{(n)}(x), \varphi_1^{(n)}(x), \cdots, \varphi_n^{(n)}(x) \right\}$$

is made, where every $\varphi_i^{(n)}(x) \in C[a,b]$, and $\varphi_0^{(n)}(x), \varphi_1^{(n)}(x), \cdots, \varphi_n^{(n)}(x)$ are linearly independent, and they meet Kronecker condition:

$$\varphi_i^{(n)}(x)\left(x_j^{(n)}\right) = \delta_{ij}, \quad i, j = 0,1,\cdots,n.$$

In $f_n(x) = \sum_{i=0}^{n} a_i^{(n)} \varphi_i^{(n)}(x)$, if we take $a_i^{(n)} = y_i^{(n)}$, then

$$f_n(x) = \sum_{i=0}^{n} y_i^{(n)} \varphi_i^{(n)}(x) = \sum_{i=0}^{n} f\left(x_i^{(n)}\right) \varphi_i^{(n)}(x)$$

This is an interpolation function as it meets interpolation condition:

$$\left(\forall j \in \{0,1,\cdots,n\}\right)\left( f_n\left(x_j^{(n)}\right) = f\left(x_j^{(n)}\right)\right).$$

Espically we take $\varphi_i^{(n)}(x) = A_i^{(n)}(x), i = 0,1,\cdots,n$, where the definition of $A_i^{(n)}(x), i = 0,1,\cdots,n$ has been expressed in (2.1). Then

$$f_n(x) = \sum_{i=0}^{n} A_i^{(n)}(x) y_i^{(n)}, \quad n = 1,2,3,\cdots \tag{3.2}$$

is just a sequence of piecewise interpolation functions.

All in all, because $f_n(x) \in \mathrm{span}\Phi(n)$, $f(x) \in C[a,b]$, and $\mathrm{span}\Phi(n)$ is a finite dimension normed linear subspace of $C[a,b]$, for any given $\varepsilon > 0$, there exists a $n \in \mathbb{N}_+$, such that the ment $f(x)$ in $C[a,b]$ can be approximated by using an element $f_n(x)$ in $\mathrm{span}\Phi(n)$, i.e., $\| f_n - f \| < \varepsilon$.

**Definition 3.1** The sequence of conditional methematical expectations $\left\{ E\left(\eta_n \middle| \xi_n = x\right)\right\}$ shown in theorem 2.1 is called a sequence of conditional methematical expectations generated by the continuous function $f(x)$. (see [1])

**Theorem 3.1** For any continuous function $f(x) \in C[a,b]$, but assuming $f(x)$ not being constant function, if $\left\{ E\left(\eta_n \middle| \xi_n = x\right)\right\}$ is a sequence of conditional methematical expectations generated by the continuous function $f(x)$, then by means of $\left\{ E\left(\eta_n \middle| \xi_n = x\right)\right\}$ we can construct a group of continuous func-

tions:

$$\Phi(n)=\left\{\varphi_0^{(n)}(x),\varphi_1^{(n)}(x),\cdots,\varphi_n^{(n)}(x)\right\},$$

where $\varphi_l^{(n)}(x)\in C[a,b], l=0,1,\cdots,n$, such that, by using the sequence of the groups of base functions $\left\{\Phi(n)\right\}$, the sequence of interpolation functions constructed as the following

$$f_n(x)=\sum_{l=0}^{n}\varphi_l^{(n)}(x)y_l^{(n)},\quad n=1,2,3,\cdots$$

uniformly converges to $f(x)$ in $[a,b]$.

**Proof** **Case 1**. Let $f(x)$ be a strictly monotone functions, and we may as well assume that $f(x)$ is a strictly monotonically increasing function, as when $f(x)$ is a strictly monotonically decreasing function, the proof is the same as the increasing status.

For any given $x\in[a,b]$, there must exist $i\in\{1,2,\cdots,n\}$, such that $x\in\left[x_{i-1}^{(n)},x_i^{(n)}\right]$. Then we have

$$R_n(x,y)=\bigvee_{k=0}^{n}\left(A_k^{(n)}(x)\cdot B_k^{(n)}(y)\right)$$
$$=\begin{cases}\left(A_{i-1}^{(n)}(x)\cdot B_{i-1}^{(n)}(y)\right)\vee\left(A_i^{(n)}(x)\cdot B_i^{(n)}(y)\right), & y\in\left[y_{i-2}^{(n)},y_{i+1}^{(n)}\right],\\ 0, & y\in[c,d]-\left[y_{i-2}^{(n)},y_{i+1}^{(n)}\right]\end{cases}$$

Now we consider the limits of Riemann sums corresponding to two integrals

$$\int_c^d yR_n(x,y)dy,\quad \int_c^d R_n(x,y)dy$$

in $E\left(\eta_n\middle|\xi_n=x\right)=\dfrac{\int_c^d yR_n(x,y)dy}{\int_c^d R_n(x,y)dy}$ as the following:

$$\int_c^d yR_n(x,y)dy=\lim_{\lambda(T_m)\to 0}\sum_{l=1}^{m}y_l^{(m)}R_n\left(x,y_l^{(m)}\right)\Delta y_l^{(m)},$$

$$\int_c^d R_n(x,y)dy=\lim_{\lambda(T_m)\to 0}\sum_{l=1}^{m}R_n\left(x,y_l^{(m)}\right)\Delta y_l^{(m)},$$

$$\Delta y_l^{(m)}=y_l^{(m)}-y_{l-1}^{(m)}, l=1,2,\cdots,m,$$

$$\lambda\left(T_m\right)=\max\left\{\Delta y_l^{(m)}\middle|l=1,2,\cdots,m\right\}$$

where $T_m$ represents a partition of $Y=[c,d]$ as the following:

$$c=y_0^{(m)}<y_1^{(m)}<\cdots<y_m^{(m)}=d\,,$$

and in the same time we should notice the expressions (see step 1 in theorem 2.1) as following:

$$x_l^{(m)}=a+lh(m),\quad h(m)=\frac{b-a}{m},$$

$$y_l^{(m)}=f\left(x_l^{(m)}\right),\quad l=0,1,\cdots,m$$

And we also should notice that $m$ and $n$ are different, which $n$ is a temporarily fixed subscript, but $m$ will tend to infinite.

Because $f(x)$ is continuous, clearly we know that $\lambda\left(T_m\right)\to 0\Rightarrow m\to\infty$, and so above expressions can be written as

$$\int_c^d yR_n(x,y)dy=\lim_{m\to\infty}\sum_{l=1}^{m}y_l^{(m)}R_n\left(x,y_l^{(m)}\right)\Delta y_l^{(m)},$$

$$\int_c^d R_n(x,y)dy=\lim_{m\to\infty}\sum_{l=1}^{m}R_n\left(x,y_l^{(m)}\right)\Delta y_l^{(m)}$$

Let $\Delta y_0^{(m)}=\Delta y_1^{(m)}$, and since $R_n(x,y)$ is a bounded function, in fact $0\le R_n(x,y)\le 1$, and $f(x)$ is also a bounded function, we have

$$\left(\exists M(x)>0\right)\left(\forall m\in\mathbb{N}\right)\left(\left|y_0^{(m)}R_n\left(x,y_0^{(m)}\right)\right|<M(x)\right).$$

So we easily know that

$$\lim_{m\to\infty}\left[y_0^{(m)}R_n\left(x,y_0^{(m)}\right)\Delta y_0^{(m)}\right]=0=\lim_{m\to\infty}\left[R_n\left(x,y_0^{(m)}\right)\Delta y_0^{(m)}\right].$$

Then we have

$$\begin{aligned}&\lim_{m\to\infty}\sum_{l=1}^{m}y_l^{(m)}R_n\left(x,y_l^{(m)}\right)\Delta y_l^{(m)}\\&=\lim_{m\to\infty}\left[y_0^{(m)}R_n\left(x,y_0^{(m)}\right)\Delta y_0^{(m)}+\sum_{l=1}^{m}y_l^{(m)}R_n\left(x,y_l^{(m)}\right)\Delta y_l^{(m)}\right]\\&=\lim_{m\to\infty}\sum_{l=0}^{m}y_l^{(m)}R_n\left(x,y_l^{(m)}\right)\Delta y_l^{(m)},\end{aligned}$$

$$\begin{aligned}&\lim_{m\to\infty}\sum_{l=1}^{m}R_n\left(x,y_l^{(m)}\right)\Delta y_l^{(m)}\\&=\lim_{m\to\infty}\left[R_n\left(x,y_0^{(m)}\right)\Delta y_0^{(m)}+\sum_{l=1}^{m}R_n\left(x,y_l^{(m)}\right)\Delta y_l^{(m)}\right]\\&=\lim_{m\to\infty}\sum_{l=0}^{m}R_n\left(x,y_l^{(m)}\right)\Delta y_l^{(m)}\end{aligned}$$

Therefore we get

$$\int_c^d yR_n(x,y)dy=\lim_{m\to\infty}\sum_{l=0}^{m}y_l^{(m)}R_n\left(x,y_l^{(m)}\right)\Delta y_l^{(m)},$$

$$\int_c^d R_n(x,y)dy=\lim_{m\to\infty}\sum_{l=0}^{m}R_n\left(x,y_l^{(m)}\right)\Delta y_l^{(m)}$$

Because of the following expressions:

$$\left(\forall x\in[a,b]\right)\left(\int_c^d R_n(x,y)dy>0\right),$$

$$\int_c^d R_n(x,y)dy=\lim_{m\to\infty}\sum_{l=0}^{m}R_n\left(x,y_l^{(m)}\right)\Delta y_l^{(m)}$$

we have the fact that, for any given $x\in[a,b]$, there exists $N(x)\in\mathbb{N}_+$, such that

$$\left(\forall m\in\mathbb{N}_+\right)\left(m>N(x)\Rightarrow\sum_{l=0}^{m}R_n\left(x,y_l^{(m)}\right)\Delta y_l^{(m)}>0\right)$$

Then when $m>N(x)$, we get

$$E\left(\eta_n\middle|\xi_n=x\right)=\frac{\int_c^d yR_n(x,y)dy}{\int_c^d R_n(x,y)dy}=\frac{\lim_{m\to\infty}\sum_{l=0}^{m}y_l^{(m)}R_n\left(x,y_l^{(m)}\right)\Delta y_l^{(m)}}{\lim_{m\to\infty}\sum_{l=0}^{m}R_n\left(x,y_l^{(m)}\right)\Delta y_l^{(m)}}$$

$$=\lim_{m\to\infty}\frac{\sum_{l=0}^{m}y_l^{(m)}R_n\left(x,y_i^{(m)}\right)\Delta y_l^{(m)}}{\sum_{j=0}^{m}R_n\left(x,y_j^{(m)}\right)\Delta y_j^{(m)}}=\lim_{m\to\infty}\sum_{l=0}^{m}\frac{R_n\left(x,y_l^{(m)}\right)\Delta y_l^{(m)}}{\sum_{j=0}^{m}R_n\left(x,y_j^{(m)}\right)\Delta y_j^{(m)}}\cdot y_l^{(m)}$$

If we let

$$\varphi_l^{(n)}(x) = \frac{R_n\left(x, y_l^{(m)}\right)\Delta y_l^{(m)}}{\sum_{j=0}^{m} R_n\left(x, y_j^{(m)}\right)\Delta y_j^{(m)}}, \ l = 0,1,\cdots,m ,$$

then above expression can be expressed as

$$\lim_{m\to\infty}\sum_{l=0}^{m}\varphi_l^{(n)}(x)y_l^{(m)} = E\left(\eta_n \middle| \xi_n = x\right). \tag{3.3}$$

Write $f_{nm}(x) = \sum_{l=0}^{m}\varphi_l^{(n)}(x)y_l^{(m)}$ , and we get a sequence with double subscripts of continuous functions as being $\left\{ f_{nm}(x)\right\}_{n,m=1}^{\infty}$ . If we put

$$\Phi_n(m) = \left\{\varphi_0^{(n)}(x), \varphi_1^{(n)}(x), \cdots, \varphi_m^{(n)}(x)\right\}$$

Then $f_{nm}(x) = \sum_{l=0}^{m}\varphi_l^{(n)}(x)y_l^{(m)}$ is just an interpolation function with base function group $\Phi_n(m)$ .

Next we especially we take $m = n$ , i.e., if we only use the diagonal elements in $\left\{ f_{nm}(x)\right\}_{n,m=1}^{\infty}$ , then we get a sequence only with single subscripts of continuous functions as being $\left\{ f_n(x)\right\}_{n=1}^{\infty}$ , where

$$f_n(x) = f_{nn}(x) = \sum_{l=0}^{n}\varphi_l^{(n)}(x)y_l^{(n)},$$
$$n = 1,2,3,\cdots$$

Let

$$\Phi(n) = \Phi_n(n) = \left\{\varphi_0^{(n)}(x), \varphi_1^{(n)}(x), \cdots, \varphi_n^{(n)}(x)\right\},$$

and then $f_n(x) = \sum_{l=0}^{n}\varphi_l^{(n)}(x)y_l^{(n)}$ is an interpolation function based on base function group $\Phi(n)$ .

For proving that $\left\{ f_n(x)\right\}_{n=1}^{\infty}$ can uniformly converge to $f(x)$ in $X = [a,b]$ . we firsily consider $n+1$ unitary functions $R_n\left(x, y_l^{(n)}\right), l = 0,1,\cdots,n$ , with respect to $x$ . In fact, since $B_i^{(n)}\left(y_l^{(n)}\right)$ are of Kronecker character:

$$B_i^{(n)}\left(y_l^{(n)}\right) = \delta_{il} = \begin{cases} 1, & i = l, \\ 0, & i \neq l, \end{cases}$$
$$i, l = 0,1,\cdots,n$$

We easily learn the following expression:

$$R_n\left(x, y_l^{(n)}\right) = \bigvee_{i=0}^{n}\left(A_i^{(n)}(x)\cdot B_i^{(n)}\left(y_l^{(n)}\right)\right) = A_l^{(n)}(x) .$$

And then we get thye following expressions:

$$\sum_{l=0}^{n} y_l R_n\left(x, y_l^{(n)}\right)\Delta y_l^{(n)} = \sum_{l=0}^{n} y_l^{(n)} A_l^{(n)}(x)\Delta y_l^{(n)},$$

$$\sum_{l=0}^{n} R_n\left(x, y_l^{(n)}\right)\Delta y_l^{(n)} = \sum_{l=0}^{n} A_l^{(n)}(x)\Delta y_l^{(n)},$$

$$\varphi_l^{(n)}(x) = \frac{R_n\left(x, y_l^{(n)}\right)\Delta y_l^{(n)}}{\sum_{j=0}^{n} R_n\left(x, y_j^{(n)}\right)\Delta y_j^{(n)}} = \frac{A_l^{(n)}(x)\Delta y_l^{(n)}}{\sum_{j=0}^{n} A_j^{(n)}(x)\Delta y_j^{(n)}},$$
$$l = 0,1,\cdots,n,$$

By these expressions we have

$$f_n(x)=\sum_{l=0}^{n}\frac{A_l^{(n)}(x)\Delta y_l^{(n)}}{\sum_{j=0}^{n}A_j^{(n)}(x)\Delta y_j^{(n)}}\cdot y_l^{(n)}=\sum_{l=0}^{n}\varphi_l^{(n)}(x)y_l^{(n)} . \tag{3.4}$$

Clearly the function group $\Phi(n)=\left\{\varphi_l^{(n)}(x)\middle| l=0,1,\cdots,n\right\}$ is linearly independent and meets the normalizing condition:

$$\left(\forall x\in[a,b]\right)\left(\sum_{l=0}^{n}\varphi_l^{(n)}(x)\equiv 1\right).$$

And the group is a dimension normed linear subspace of $C[a,b]$. If $\Phi(n)$ is regarded as a base function group, then

$$f_n(x)=\sum_{l=0}^{n}\varphi_l^{(n)}(x)y_l^{(n)} \tag{3.5}$$

is just a piecewise interpolation function.

At last, we prove that the sequence $\left\{f_n(x)\right\}_{n=1}^{\infty}$ uniformly converges to $f(x)$ in $[a,b]$.

In fact, for any given $x\in[a,b]$, clearly $\left(\exists i\in\{1,2,\cdots,n\}\right)\left(x\in\left[x_{i-1}^{(n)},x_i^{(n)}\right]\right)$, and then

$$f_n(x)=\sum_{l=0}^{n}\varphi_l^{(n)}(x)y_l^{(n)}=\varphi_{i-1}^{(n)}(x)y_{i-1}^{(n)}+\varphi_i^{(n)}(x)y_i^{(n)} .$$

Because $f(x)$ is continuous, there exist two points $\xi_i^{(n)},\eta_i^{(n)}\in\left[x_{i-1}^{(n)},x_i^{(n)}\right]$, such that

$$f\left(\xi_i^{(n)}\right)=\min_{x\in\left[x_{i-1}^{(n)},x_i^{(n)}\right]}f(x),\quad f\left(\eta_i^{(n)}\right)=\max_{x\in\left[x_{i-1}^{(n)},x_i^{(n)}\right]}f(x)$$

Clearly the following expressions are true:

$$\left(\forall l\in\{0,1,\cdots,n\}\right)\left(\varphi_l^{(n)}\left([a,b]\right)=[0,1]\right),$$

$$\left(\forall x\in\left[x_{i-1}^{(n)},x_i^{(n)}\right]\right)\left(\varphi_{i-1}^{(n)}(x)+\varphi_i^{(n)}(x)=1\right).$$

So we have the result: for any $x\in\left[x_{i-1}^{(n)},x_i^{(n)}\right]$, we have

$$f\left(\xi_i^{(n)}\right)\le\varphi_{i-1}^{(n)}(x)y_{i-1}^{(n)}+\varphi_i^{(n)}y_i^{(n)}\le f\left(\eta_i^{(n)}\right).$$

By noticing the facts $y_{i-1}^{(n)}=f\left(x_{i-1}^{(n)}\right)$ and $y_i^{(n)}=f\left(x_i^{(n)}\right)$, for any point $x\in\left[x_{i-1}^{(n)},x_i^{(n)}\right]$, we have the following inequation:

$$\begin{aligned}\left|f(x)-f_n(x)\right|&=\left|f(x)-\left(\varphi_{i-1}^{(n)}(x)y_{i-1}^{(n)}+\varphi_i^{(n)}y_i^{(n)}\right)\right|\\&\le\left|f\left(\eta_i^{(n)}\right)-f\left(\xi_i^{(n)}\right)\right|\end{aligned}$$

Also for $f(x)$ being uniformly continuous in $[a,b]$, for any given $\varepsilon>0$, there exists $\delta>0$, such that

$$\left(\forall u,v\in[a,b]\right)\left(|u-v|<\delta\Rightarrow\left|f(u)-f(v)\right|<\varepsilon\right).$$

Now we suppose that $h(n)=\dfrac{b-a}{n}<\delta$, then $\left|f\left(\eta_i^{(n)}\right)-f\left(\xi_i^{(n)}\right)\right|<\varepsilon$, and we have

$$\left(\forall x\in\left[x_{i-1}^{(n)},x_i^{(n)}\right]\right)\left(\left|f(x)-f_n(x)\right|<\varepsilon\right).$$

As $h(n)\to 0\Leftrightarrow n\to\infty$, so

$$\left(\exists N\in\mathbb{N}_+\right)\left(\forall n\in\mathbb{N}_+\right)\left(n>N\Rightarrow h(n)<\delta\right);$$

hence

$$\left(\forall n\in\mathbb{N}_+\right)\left(n>N\Rightarrow\left(\forall x\in[a,b]\right)\left(\left|f(x)-f_n(x)\right|<\varepsilon\right)\right).$$

This means that $\left\{f_n(x)\right\}_{n=1}^{\infty}$ uniformly converges to $f(x)$ in $X=[a,b]$.

**Case 2**. Assume $f(x)$ be not strictly monotonic. Similar to the method in theorem 2.1, we can have

$$R_n(x,y) = \bigvee_{j=0}^{n} \left( A_{k_j}^{(n)}(x) \cdot B_{k_j}^{(n)}(y) \right).$$

Based on them, we get a sequence of interpolation functions:

$$f_n(x) = \sum_{l=0}^{n} \frac{A_{k_l}^{(n)}(x)\Delta y_{k_l}^{(n)}}{\sum_{s=0}^{n} A_{k_s}^{(n)}(x)\Delta y_{k_s}^{(n)}} \cdot y_{k_l}^{(n)}$$

$$= \sum_{l=0}^{n} \varphi_l^{(n)}(x) y_{k_l}^{(n)},$$

where

$$\varphi_l^{(n)}(x) = \frac{A_{k_l}^{(n)}(x)\Delta y_{k_l}^{(n)}}{\sum_{s=0}^{n} A_{k_s}^{(n)}(x)\Delta y_{k_s}^{(n)}},$$

$$l = 0,1,\cdots,n$$

Then by the same way we know that $\{ f_n(x) \}_{n=1}^{\infty}$ uniformly converges to $f(x)$ in $X = [a,b]$. The proof of the theorem has been finished. □

## 4. Quantum Mechanics Representation of Classic Mechanics

As we all know, classic mechanics is the scope of macroscopical physics in which Newtonian mechanics is its main part. Classic mechanics is very different with microphysics, especially with quantum mechanics (see [2][3]). For example, the motions of microscopic particles have wave-particle duality; but the motion of mass points in macroscopical physics only has mass point characters and no wave natures; in other words, there is no wave-mass-point duality in macroscopical physics. However, there has ever exisited a correspondence principle: considering a kind of motion state in quantum physics, when quantum number $n \to \infty$, the limit situation of the motion state in quantum physics must become a kind of motion state in macroscopical physics. In other words, the limit situation of the motion law in quantum physics is just some motion law in macroscopical physics.

Generally, Bohr suggested a generalized correspondence principle: the limit situation of any new theory must be in line with some old theory.

It is worth noting that above correspondence principle or generalized correspondence principle is all of unipolarity: the limit situation of the motion law in quantum physics is just some motion law in macroscopical physics, but the converse principle is clearly meaningless.

However, we can consider an important problem: must any one of motion states in macroscopical physics be the limit situation of some the motion states in quantum physics?

Apparently, this problem has not been observed, and of course there is no any answer. For example, we consider the well-known projectile motion. As we all know, a projectile motion can be expressed by the equation of locus of the projectile motion as following:

$$y(x) = x\tan\alpha - \frac{g}{2v_0^2\cos^2\alpha} x^2,$$

$$x \in [0, d_0], \quad d_0 = \frac{2v_0^2}{g}\sin 2\alpha$$

where $\alpha \in \left(0, \frac{\pi}{2}\right)$ is a mass ejection angle, $d_0 \in (0,+\infty)$ is the maximum range of fire, and the initial velocity is repressed by $v_0 \in (0,+\infty)$; here the air friction is omitted.

Clearly $y(x) \in C[0,d_0]$, i.e., a projectile motion can be described by a unary continuous function. For this continuous function, can we find some microscopic particles such that the limit of the group behavior of these microscopic particles is just this continuous function $y(x)$ when quantum number $n \to \infty$ ?

In this paper, we will give a positive answer for this problem. It is easy to understand that almost all laws of classic mechanics are described by continuous functions. So we can generalize above problem as such problem: for any a continuous function $f$ , unary continuous function or multivariate continuous function, which should describe some motion law of some mass point in microscopically physics, can we find some microscopic particles such that the limit of the group behavior of these microscopic particles is just this continuous function $f$ when quantum number $n \to \infty$ ?

Now we start to try to solve out the problem.

Firstly, we consider the case of unary continuous functions. For any a unary continuous function $f(x) \in C[a,b]$, we can use a linear transformation as the following:

$$u:[a,b] \to [0,1], x \mapsto u = u(x) = \frac{x-a}{b-a}$$

to redefine the continuous function $f(x)$ on the closed interval $[0,1]$, i.e.,

$$g(u) = f\left((b-a)u + a\right) \in C[0,1],$$

Therefore, without loss of generality, we can only consider such continuous functions as being $f(x) \in C[0,1]$. However, we do not consider constant functions because constant functions are almost meaningless in physics.

**Theorem 4.1** Given arbitrarily a non-constant function $f(x) \in C[0,1]$, there must exist some microscopic particles such that the limit of the group behavior of these microscopic particles is just this continuous function $f(x)$ when the quantum number $n \to \infty$ .

**Proof**. **Step 1**. We consider the wave function of a microscopic particle in infinite deep square potential well.

As a matter of fact, we take a particle $M$ with quality $m$ , and $M$ moves along $Ox$ axis and is of determined momentum $p = mv_x$ and determined energy $E = \frac{1}{2}mv_x^2 = \frac{p^2}{2m}$ where $v_x$ is the velocity of $M$ moving along $Ox$ . We take a special infinite deep square potential well as follows (see Fig. 4.1):

$$V(x) = \begin{cases} 0, & x \in [0,1], \\ +\infty, & x \in (-\infty, 0) \cup (1, +\infty) \end{cases}$$

The particle $M$ is complete free inside the potential well; only at two endpoints $x = 0, x = 1$, there are infinite forces to impose restrictions on $M$ not to escape.

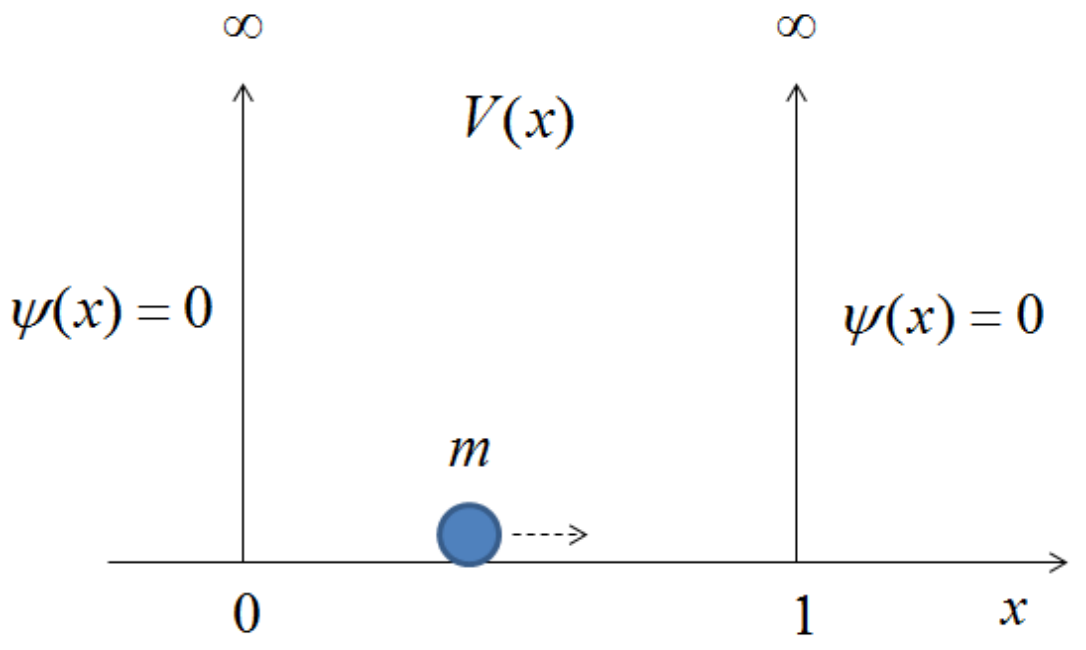

Fig. 4.1. Particle movement in the infinite deep square potential well

At outside of the potential well, i.e., $x \in (-\infty, 0) \cup (1, +\infty)$ , we now notice the steady state Schrodinger Equation as the following:

$$\left[ -\frac{\hbar^2}{2m} \frac{\partial^2}{\partial x^2} + V(x) \right] \psi(x) = E\psi(x)$$

It is easy to know that $\psi(x) = 0$ ; so the probability of finding the particle in the interval

$(-\infty,0)\bigcup(1,+\infty)$ is zero. However inside the potential well, i.e., $x\in[0,1]$, we have $V(x)=0$; then the Schrodinger Equation turn into the following form:

$$\frac{d^2\psi(x)}{dx^2}=-\left(\frac{\sqrt{2mE}}{\hbar}\right)^2\psi(x).$$

Let $k=\frac{\sqrt{2mE}}{\hbar}$; then we have the following form:

$$\frac{d^2\psi(x)}{dx^2}+k^2\psi(x)=0,$$

which is the equation of motion of a simple harmonic oscillation, and its general solution is as following:

$$\psi(x)=A\sin kx+B\cos kx,$$

where $A,B$ are two arbitrary constants that can be determined by some boundary conditions.

Then, what are the boundary conditions? In fact, in quantum mechanics, the solution of a Schrodinger Equation in three-dimension space, i.e. the wave function $\Psi(x,y,z,t)$ should satisfy the following established standard conditions:

i) $\int_{\mathbb{R}^3}|\Psi|^2dxdydz=1$;

ii) $\Psi$ and its three partial derivatives $\frac{\partial\Psi}{\partial x},\frac{\partial\Psi}{\partial y},\frac{\partial\Psi}{\partial z}$ are continuous everywhere;

iii) $\Psi$ is a single-valued function about coordinates.

By means of above conditions, when the potential function approaches infinite, based on the continuity of $\psi(x)$, we can get the result as being $\psi(0)=\psi(1)=0$, which can make the solution be continuous at both inside and outside of the potential well. Because of the following expression:

$$0=\psi(0)=A\sin 0+B\cos 0=B,$$

we get $B=0$; thus we have the following equation:

$$\psi(x)=A\sin kx$$

And then we take notice of the equation: $0=\psi(1)=A\sin k$, if $A=0$, then $\psi(x)\equiv 0$ which is a trivial solution and cannot be normalized. Thus we only get the result: $\sin k=0$, and we know the following fact:

$$k=0,\pm\pi,\pm2\pi,\pm3\pi,\cdots$$

Clearly that $k=0$ is meaningless, for this can also make that $\psi(x)\equiv 0$. Besides, $k$ with negative values cannot generate any new solutions because of the fact that $\sin(-\theta)=-\sin\theta$ and we can make the minus signs enter into the coefficient $A$. Therefore, we have the result:

$$k=k_n=n\pi,\quad n=1,2,3,\cdots$$

We should notice a fact that, the boundary condition at $x=1$ is not used to determine the coefficient $A$, but to determine the energy $E$ because of the expression: $k=\frac{\sqrt{2mE}}{\hbar}$, i.e.,

$$E=E_n=\frac{\hbar^2k_n^2}{2m}=\frac{n^2\pi^2\hbar^2}{2m},\quad n=1,2,3,\cdots \tag{4.1}$$

It is well-known that $E_1=\frac{\pi^2\hbar^2}{2m}$ is ground state, and others are follows:

$$E_2=4E_1,E_3=9E_1,E_4=16E_1,\cdots$$

which means that the energy of a particle can only take discrete values; in other words, the energy of a particle is quantized. And positive integer $n$ is called the quantum number of the energy of a particle. So we can learn that the quantization of the energy of a particle is very natural in quantum mechanics.

Thus the solution of the Schrodinger Equation can be expressed by the quantum number as the following:

$$\psi_n(x) = A\sin(n\pi x), \quad n = 1,2,3,\cdots \tag{4.2}$$

In order to determine the coefficient $A$, we can use the normalization condition $\int_0^1 |\psi_n(x)|^2 dx = 1$ to get $A = \sqrt{2}$. Then we get the solution of the Schrodinger Equation inside the potential well as the following:

$$\begin{aligned} &\psi_n(x) = \sqrt{2}\sin(n\pi x), \\ &x \in [0,1], \quad n = 1,2,3,\cdots \end{aligned} \tag{4.3}$$

Let

$$\alpha_n(x) = \sin(n\pi x), \quad x \in [0,1], \quad n = 1,2,3,\cdots,$$

and we have the following form:

$$\begin{aligned} &\psi_n(x) = \sqrt{2}\alpha_n(x), \\ &x \in [0,1], \quad n = 1,2,3,\cdots \end{aligned} \tag{4.4}$$

The function $\alpha_n(x)$ is called **essence wave function** of the wave function $\psi_n(x)$.

**Step 2.** Based on an important fact that will be described as follows, we should consider to weaken three standard conditions about the wave function $\Psi(x,y,z,t)$ mentioned above.

As a matter of fact, we can see that the derived function $\frac{\partial\psi_n(x)}{\partial x}$ of the wave function $\psi_n(x) = \sqrt{2}\sin(n\pi x)$ is not continuous at $x = 0,1$. For this we only notice the following implication is true:

$$\begin{aligned} &\frac{\partial\psi_n(x)}{\partial x} = \frac{\sqrt{2}}{n\pi}\cos(n\pi x) \Rightarrow \\ &\left( \frac{\partial\psi_n(0)}{\partial x} = \frac{\sqrt{2}}{n\pi} \neq 0, \quad \frac{\partial\psi_n(1)}{\partial x} = \frac{\sqrt{2}}{n\pi}\cos(n\pi) \neq 0 \right). \end{aligned}$$

It is well known that the movement of a particle in the infinite potential well is a typical example in quantum mechanics. However, as we have learned above, its wave function $\Psi$ and its three partial derivatives $\frac{\partial\Psi}{\partial x}, \frac{\partial\Psi}{\partial y}, \frac{\partial\Psi}{\partial z}$ are not continuous at everywhere (of course, in above case, there is only one partial derivative $\frac{\partial\Psi}{\partial x}$, in fact $\frac{\partial\Psi}{\partial x} = \frac{\mathrm{d}\Psi}{\mathrm{d}x}$ here).

We should not forget a fact that wave function $\Psi$ does not represent a physical wave but only a mathematical wave; in other words, $|\Psi|^2$ is a probability density function where it should be normalized.

We also know such a fact that, in probability theory, any probability density function is not required to be continuous at everywhere but only required to be almost everywhere continuous. Thus, we have enough reason to revise the three standard conditions which the wave function $\Psi(x,y,z,t)$ should satisfy mentioned above to be as the following:

(i) $\int_{\mathbb{R}^3} |\Psi|^2 dxdydz = 1$;

(ii) $\Psi$ and its three partial derivatives $\frac{\partial\Psi}{\partial x}, \frac{\partial\Psi}{\partial y}, \frac{\partial\Psi}{\partial z}$ cannot be continuous only at finite points (clearly the requirement is a little stronger than almost everywhere continuous);

(iii) $\Psi$ is a single-valued function about coordinates.

Moreover, by the viewpoint of Von Neumann, wave function $\Psi$ is defined in a Hilbert space $\mathcal{L}^2(\mathbb{R}^3)$, where the operations in quantum mechanics (momentum, work, and so on) are inner product operations, which may be enlightened by $\int_{\mathbb{R}^3} |\Psi|^2 dxdydz = 1$ and form a mathematical formalization structure. We all know the fact that, in a Hilbert space $\mathcal{L}^2(\mathbb{R}^3)$, we have no need to require wave function $\Psi$ to be continuous at everywhere but almost everywhere continuous to be enough.

**Step 3**. We continue to consider the wave function of the particle in the one dimension infinite deep potential well. We have known its general solution being as

$$\psi(x) = A\sin kx + B\cos kx,$$

where $A, B$ are arbitrary constants which can be determined by the boundary conditions. This time, we suppose $\frac{\partial\psi(x)}{\partial x}$ be continuous at the boundary points $x = 0,1$. We take notice of the following implication:

$$\frac{\partial\psi(x)}{\partial x} = \frac{A}{k}\cos kx - \frac{B}{k}\sin kx$$

$$\Rightarrow 0 = \frac{\partial\psi(0)}{\partial x} = \frac{A}{k} \Rightarrow A = 0$$

Then we get the following result:

$$\psi(x) = B\cos kx.$$

And then we pay attention to the equation $\frac{\partial\psi(x)}{\partial x} = -\frac{B}{k}\sin kx$, so that

$$0 = \frac{\partial\psi(1)}{\partial x} = -\frac{B}{k}\sin k.$$

Because $\frac{B}{k} \neq 0$, we solve out the values of $k$ as follows:

$$k = k_n = n\pi, \quad n = 1,2,3,\cdots$$

Very similar to the method in Step 1, we have the expression of $E$ again as follows:

$$E = E_n = \frac{\hbar^2 k_n^2}{2m} = \frac{n^2\pi^2\hbar^2}{2m}, \quad n = 1,2,3,\cdots$$

So the solution of the Schrodinger Equation can be expressed by means of quantum numbers as the following:

$$\varphi_n(x) = B\cos(n\pi x), \quad n = 1,2,3,\cdots \tag{4.5}$$

Again by using the normalization condition, we can get that $B = \sqrt{2}$. Thus another solution of the Schrodinger Equation in the potential well is as the following:

$$\begin{gathered} \varphi_n(x) = \sqrt{2}\cos(n\pi x), \\ x \in [0,1], \quad n = 1,2,3,\cdots \end{gathered} \tag{4.6}$$

Let $\beta_n(x) = \cos(n\pi x), x \in [0,1], n = 1,2,3,\cdots$, and we have

$$\begin{gathered} \varphi_n(x) = \sqrt{2}\beta_n(x), \\ x \in [0,1], \quad n = 1,2,3,\cdots \end{gathered} \tag{4.7}$$

The function $\beta_n(x)$ is also called **essence wave function** of the wave function $\varphi_n(x)$.

It is interesting to note that the wave function $\varphi_n(x) = \sqrt{2}\cos(n\pi x)$ is not continuous at boundary points $x = 0,1$ this time. Besides, since

$$\psi_n\left(x + \frac{1}{2n}\right) = \sqrt{2}\sin\left[n\pi\left(x + \frac{1}{2n}\right)\right]$$

$$= \sqrt{2}\sin\left(n\pi x + \frac{\pi}{2}\right) = \sqrt{2}\cos(n\pi x) = \varphi_n(x)$$

when the quantum number $n$ is very large, the two wave functions $\psi_n(x)$ and $\varphi_n(x)$ are almost no different; in other words, $\varphi_n(x)$ is just the situation that $\psi_n(x)$ translates a $\frac{\pi}{2}$ phase position to the right side.

For visualization, the function $\psi_n(x)$ can be vividly called Adam wave function and $\varphi_n(x)$ be called Eve wave function. In fact, we care more about the function family of essence wave functions of Adam and Eve wave functions, denoted by $\{\alpha_n(x), \beta_n(x)\}_{n=1}^{\infty}$, and we can call $\alpha_n(x)$ to be Adam essence wave function and $\beta_n(x)$ to be Eve essence wave function. Clearly $\alpha_n(x)$ and $\beta_n(x)$ are defined on the unit interval $X = [0,1]$.

**Step 4**. Supplementary instruction for the revision of the three standard requirements on the wave function $\Psi$

It is well known that, in physics, harmonic oscillation is often described by complex exponential form; for example, the two wave functions that we just get can be described as the following:

$$\begin{aligned}\Psi(x) &= \sqrt{2}e^{i(n\pi x)} \\ &= \sqrt{2}\cos(n\pi x) + i\sqrt{2}\sin(n\pi x) \\ &= \varphi_n(x) + i\psi_n(x)\end{aligned} \tag{4.8}$$

In classic physics, this kind of expression is said to be more convenient for operation but with no more physical significance. However, here we can find its physical significance of the complex variables function $\Psi(x) = \sqrt{2}e^{in\pi x}$ coming from quantum mechanics. As its real part of the $\Psi(x) = \sqrt{2}e^{in\pi x}$, Eve wave function as being $\varphi_n(x) = \sqrt{2}\cos(n\pi x)$ is determined by the second boundary condition; and its imaginary part, Adam wave function as being $\psi_n(x) = \sqrt{2}\sin(n\pi x)$ is determined by the first boundary condition. These mean that the two boundary conditions are all useful and we cannot give up any one of them. Therefore, the revision of the three standard requirements is quite reasonable.

**Step 5**. The extension of the domain of definition of the wave functions

For any a finite closed interval $[a,b]$, by means of the linear transformation as follows:

$$t = (b-a)x + a\,,$$

the essence wave function family $\{\sin(n\pi x), \cos(n\pi x)\}_{n=1}^{\infty}$ defined on the interval $[0,1]$ can be extended to the closed interval $[a,b]$; we rewrite the variable $t$ to be $x$, and we have the following form:

$$\left\{\sin\frac{n\pi(x-a)}{b-a}, \cos\frac{n\pi(x-a)}{b-a}\right\}_{n=1}^{\infty}, \quad x \in [a,b]$$

We can easily know that the mapping as follows

$$u : [a,b] \to [0,1], \quad x \mapsto u(x) = \frac{x-a}{b-a}$$

is a topological homeomorphism from $[a,b]$ to $[0,1]$. This means that the essence wave function family $\{\sin(n\pi x), \cos(n\pi x)\}_{n=1}^{\infty}$ and the family of essence wave functions

$$\left\{\sin\frac{n\pi(x-a)}{b-a}, \cos\frac{n\pi(x-a)}{b-a}\right\}_{n=1}^{\infty}$$

is not essentially different; so they can be regarded the same.

It is worth noting that, for Adam wave function, in the infinite deep square potential well, it should be written as the following complete form:

$$\Psi(x,t) = \psi_n(x)e^{-\frac{i}{\hbar}E_n t} = \sqrt{2}\sin(n\pi x)e^{-\frac{i}{\hbar}E_n t}, \quad x \in [0,1] \tag{4.9}$$

where we only write out the expression just as being $x \in [0,1]$. And for Eve wave function, in the infinite deep square potential well, it should be written as the following complete form:

$$\Psi(x,t)=\varphi_n(x)e^{-\frac{i}{\hbar}E_nt}=\sqrt{2}\cos(n\pi x)e^{-\frac{i}{\hbar}E_nt}, \quad x\in[0,1] \tag{4.10}$$

where we also only write out the expression just as being $x\in[0,1]$.

Based on the statistical interpretation of wave functions, $|\Psi(x,t)|^2$ should be a kind of probability density function. Then from Equations (4.9) and (4.10), we can learn that $2\sin^2(n\pi x)$ is a probability density function and $2\cos^2(n\pi x)$ is a probability density function too. We have enough reason to call $\sin^2(n\pi x)$ and $\cos^2(n\pi x)$ essence probability density function of the probability density functions. So we get the essence probability density function family of Adam wave functions and Eve functions as the following:

$$\left\{\sin^2(n\pi x),\cos^2(n\pi x)\right\}_{n=1}^{\infty}. \tag{4.11}$$

It is easy to see that $\left\{\sin^2(n\pi x),\cos^2(n\pi x)\right\}_{n=1}^{\infty}$ is of two-phase normalization property:

$$\sin^2(n\pi x)+\cos^2(n\pi x)=1.$$

**Step 6**. The construction of the sequence of two-dimension probability density functions

Given arbitrarily a continuous function $f\in C[0,1]$, clearly $f\left([0,1]\right)$ is a closed interval, denoted by $Y=[c,d]=f\left([0,1]\right)$. Let

$$X(n)=\left\{x\in[0,1]\middle|\sin(n\pi x)=0,\cos(n\pi x)=0\right\}.$$

And we can easily know that $\mathrm{card}\left(X(n)\right)=2n+1$. Hence we have the following expression:

$$X(n)=\left\{x_i^{(n)}\middle|i=0,1,2,\cdots,2n\right\},$$

where $x_i^{(n)}=\dfrac{i}{2n}$, $i=0,1,2,\cdots,2n$. This expression means that the closed interval $X=[0,1]$ are equidistantly partitioned as the following:

$$\Delta x_i^{(n)}=x_i^{(n)}-x_{i-1}^{(n)}=\frac{1}{2n},\quad i=1,2,\cdots,2n.$$

And we let

$$Y(n)=\left\{y_i^{(n)}=f\left(x_i^{(n)}\right)\middle|i=0,1,2,\cdots,2n\right\}.$$

For convenience, let $m=2n$; but be careful, here $m$ means subscript but not the quality of some particle. We are going to discuss our problem from the following two cases.

**Case 1**. Suppose $f(x)$ is a strict monotonous function. It assumes that $f(x)$ be a strict monotonous rising function because its proof is not of essence difference when $f(x)$ is a strict monotonous declining function. Therefore, we have the following partition:

$$c=y_0^{(n)}<y_1^{(n)}<\cdots<y_{m-1}^{(n)}<y_m^{(n)}=d.$$

Then we consider the particle wave functions defined in the following subintervals one by one:

$$\left[y_0^{(n)},y_1^{(n)}\right],\left[y_1^{(n)},y_2^{(n)}\right],\cdots,\left[y_{m-2}^{(n)},y_{m-1}^{(n)}\right],\left[y_{m-1}^{(n)},y_m^{(n)}\right].$$

Firstly we treat with it in the closed interval $\left[y_0^{(n)},y_1^{(n)}\right]$. And we consider the movement of a particle in the infinite deep square potential well that the closed interval $\left[0,2\left(y_1^{(n)}-y_0^{(n)}\right)\right]$ is just the bottom margin of the potential well. The particle is denoted by $M_1^{(n)}$ which can be regarded as a descendant particle generated by the Adam wave function and Eve wave function of the original particle $M$ in the case of energy level being $n$. The descendant particle $M_1^{(n)}$ moves along $Oy$ axis with determined quality $m_1^{(n)}$ and determined momentum $p_1^{(n)}=m_1^{(n)}v_y^{(n,1)}$ and determined energy

$$E_1^{(n)} = \frac{1}{2} m_1^{(n)} \left( v_y^{(n,1)} \right)^2 = \frac{\left( p_1^{(n)} \right)^2}{2m_1^{(n)}},$$

where $v_y^{(n,1)}$ is the velocity of movement of $M_1^{(n)}$ along $Oy$ axis. By means of the continuity of the wave function, it is easy to get the solution of the wave function in $\left[0, 2\left(y_1^{(n)} - y_0^{(n)}\right)\right]$ as following:

$$\psi_p^{(n,1)}(y) = \sqrt{\frac{2}{2\left(y_1^{(n)} - y_0^{(n)}\right)}} \sin \frac{p\pi}{2\left(y_1^{(n)} - y_0^{(n)}\right)} y, \quad p = 1, 2, 3, \cdots \tag{4.12}$$

Then again, by means of the continuity of the derived function of the wave function, we can get another solution of the wave function in the closed interval $\left[0, 2\left(y_1^{(n)} - y_0^{(n)}\right)\right]$ as following:

$$\varphi_p^{(n,1)}(y) = \sqrt{\frac{2}{2\left(y_1^{(n)} - y_0^{(n)}\right)}} \cos \frac{p\pi}{2\left(y_1^{(n)} - y_0^{(n)}\right)} y, \quad p = 1, 2, 3, \cdots \tag{4.13}$$

Now we care more for the ground state of $\psi_p^{(n,1)}(y)$ and $\varphi_p^{(n,1)}(y)$, i.e., the wave functions when $p = 1$ as follows:

$$\psi_1^{(n,1)}(y) = \sqrt{\frac{2}{2\left(y_1^{(n)} - y_0^{(n)}\right)}} \sin \frac{\pi}{2\left(y_1^{(n)} - y_0^{(n)}\right)} y, \tag{4.14}$$

$$\varphi_1^{(n,1)}(y) = \sqrt{\frac{2}{2\left(y_1^{(n)} - y_0^{(n)}\right)}} \cos \frac{\pi}{2\left(y_1^{(n)} - y_0^{(n)}\right)} y, \tag{4.15}$$

We can omit the amplitudes of the wave and keep the essence wave function and do squaring operation on the essence wave functions, and get the probability essence wave functions as the following:

$$\sin^2 \frac{\pi}{2\left(y_1^{(n)} - y_0^{(n)}\right)} y, \quad \cos^2 \frac{\pi}{2\left(y_1^{(n)} - y_0^{(n)}\right)} y. \tag{4.16}$$

The graphs of the probability essence wave functions in $\left[0, y_1^{(n)} - y_0^{(n)}\right]$ are shown in Fig. 4.2.

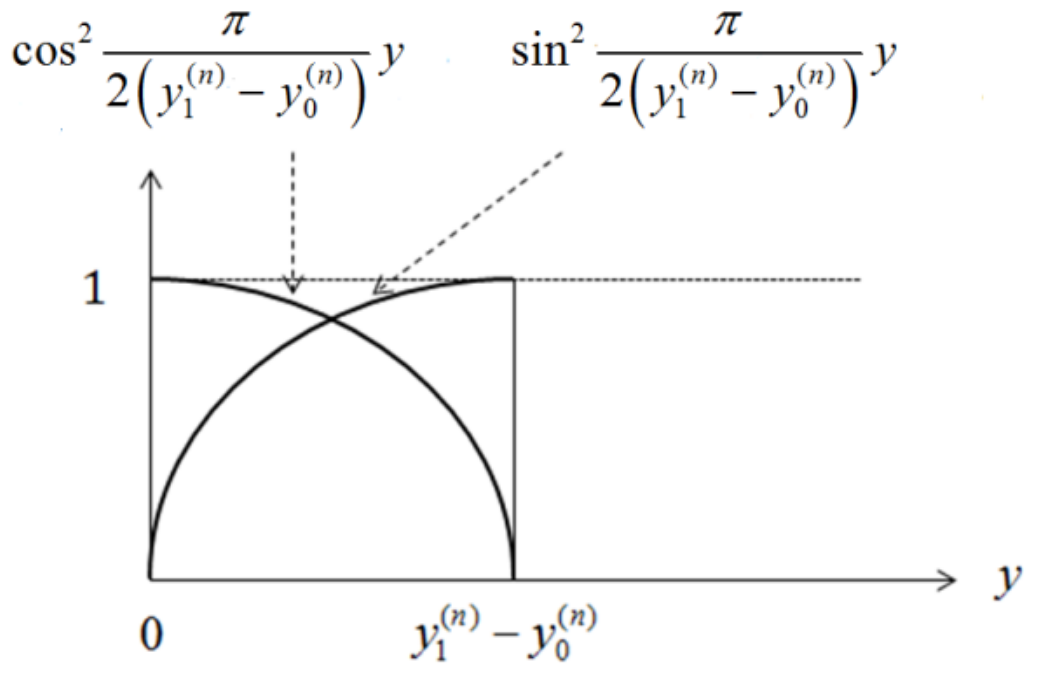

Fig. 4.2. The probability essence wave functions in $\left[0, y_1^{(n)} - y_0^{(n)}\right]$

The next, we make a coordinate translation: $t = y + y_0^{(n)}$, and then we have the following expressions:

$$\sin^2 \frac{\pi}{2\left(y_1^{(n)} - y_0^{(n)}\right)} y = \sin^2 \frac{\pi}{2\left(y_1^{(n)} - y_0^{(n)}\right)} \left(t - y_0^{(n)}\right),$$

$$\cos^2\frac{\pi}{2\left(y_1^{(n)}-y_0^{(n)}\right)}y=\cos^2\frac{\pi}{2\left(y_1^{(n)}-y_0^{(n)}\right)}\left(t-y_0^{(n)}\right)$$

Thus we transfer the probability essence wave functions defined in the closed interval $\left[0,y_1^{(n)}-y_0^{(n)}\right]$ into the probability essence wave functions in in closed interval $\left[y_0^{(n)},y_1^{(n)}\right]$. And we rewrite the variable $t$ back to $y$, and then we get the following expressions:

$$\sin^2\frac{\pi}{2\left(y_1^{(n)}-y_0^{(n)}\right)}\left(y-y_0^{(n)}\right),\quad \cos^2\frac{\pi}{2\left(y_1^{(n)}-y_0^{(n)}\right)}\left(y-y_0^{(n)}\right) \tag{4.17}$$

The graphs of the probability essence wave functions in $\left[y_0^{(n)},y_1^{(n)}\right]$ are shown in Fig. 4.3.

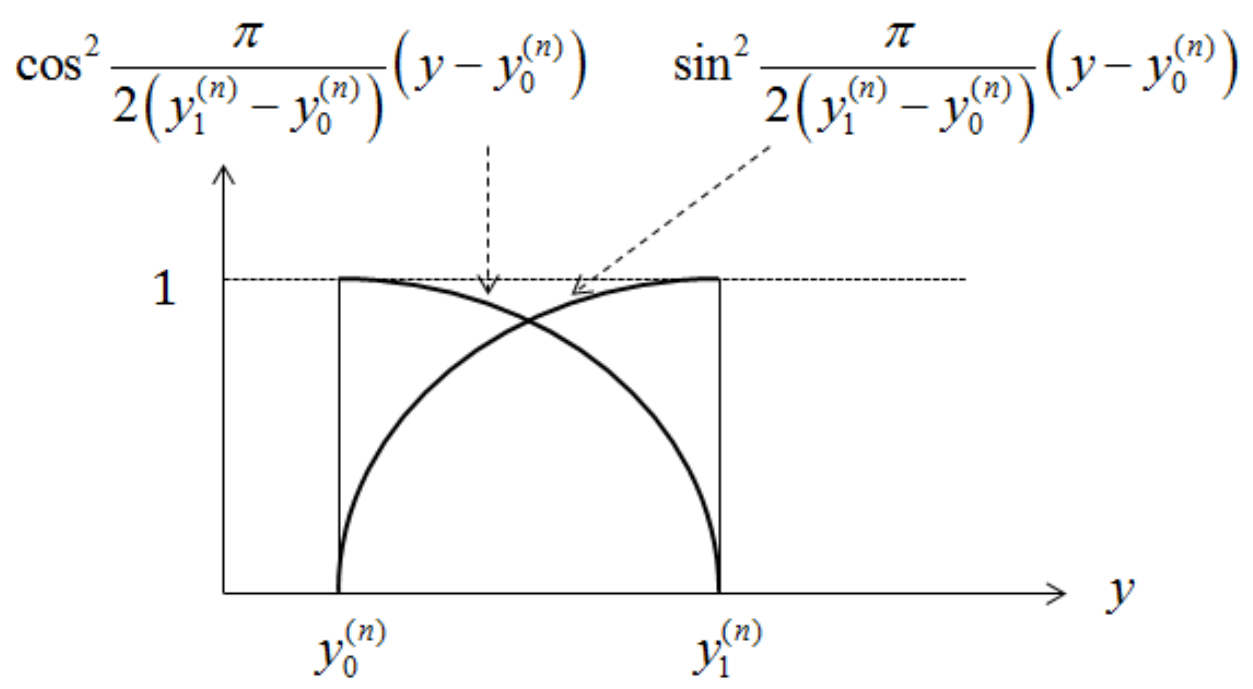

Fig. 4.3. The probability essence wave functions in $\left[y_0^{(n)},y_1^{(n)}\right]$

In that way, we can regard the following expression

$$\sin^2\frac{\pi}{2\left(y_1^{(n)}-y_0^{(n)}\right)}\left(y-y_0^{(n)}\right)$$

as Adam probability essence wave function of the movement of the descendant particle $M_1^{(n)}$ in $\left[y_0^{(n)},y_1^{(n)}\right]$, and regard the following expression

$$\cos^2\frac{\pi}{2\left(y_1^{(n)}-y_0^{(n)}\right)}\left(y-y_0^{(n)}\right)$$

as Eve probability essence wave function of the movement of the descendant particle $M_1^{(n)}$ in $\left[y_0^{(n)},y_1^{(n)}\right]$.

In the same way, we can get Adam and Eve probability essence wave functions of the movement of the descendant particles $M_2^{(n)},\cdots,M_m^{(n)}$ in the closed intervals $\left[y_1^{(n)},y_2^{(n)}\right],\cdots,\left[y_{m-1}^{(n)},y_m^{(n)}\right]$ respectively as the following:

$$\sin^2\frac{\pi}{2\left(y_2^{(n)}-y_1^{(n)}\right)}\left(y-y_1^{(n)}\right),\quad \cos^2\frac{\pi}{2\left(y_2^{(n)}-y_1^{(n)}\right)}\left(y-y_1^{(n)}\right),$$

……

$$\sin^2\frac{\pi}{2\left(y_m^{(n)}-y_{m-1}^{(n)}\right)}\left(y-y_{m-1}^{(n)}\right),\quad \cos^2\frac{\pi}{2\left(y_m^{(n)}-y_{m-1}^{(n)}\right)}\left(y-y_{m-1}^{(n)}\right)$$

We get all these graphs of the probability essence wave functions together in the following closed intervals:

$$\left[y_0^{(n)},y_1^{(n)}\right],\left[y_1^{(n)},y_2^{(n)}\right],\cdots,\left[y_{m-1}^{(n)},y_m^{(n)}\right]$$

and they are shown in Fig. 4.4.

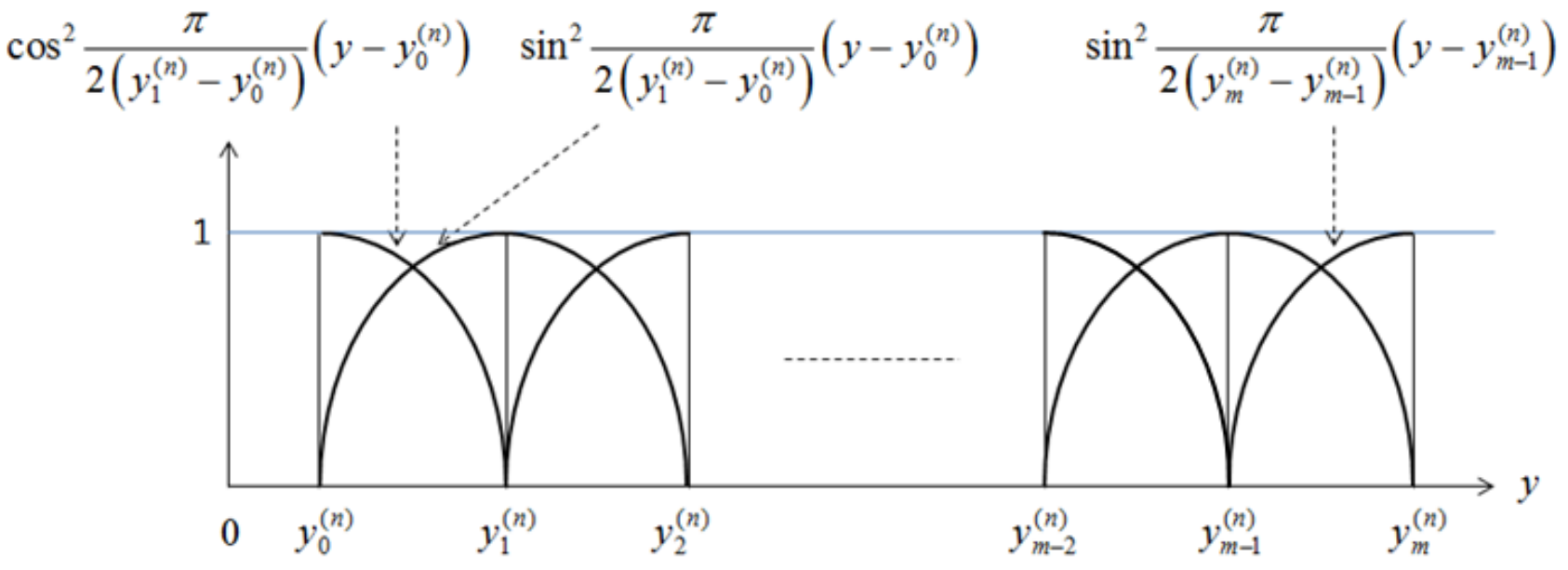

Fig. 4.4. All the probability essence wave functions in $\left[y_0^{(n)},y_1^{(n)}\right],\cdots,\left[y_{m-1}^{(n)},y_m^{(n)}\right]$

Now we need to summarize the work that we have done as follows.

When $x\in\left[x_0^{(n)},x_1^{(n)}\right]$, from the information of Adam and Eve probability essence wave functions $\sin^2(n\pi x)$ and $\cos^2(n\pi x)$ at the nodes as the following:

$$y_0^{(n)}=f\left(x_0^{(n)}\right),\quad y_1^{(n)}=f\left(x_1^{(n)}\right),$$

we get Adam and Eve probability essence wave functions in $\left[y_0^{(n)},y_1^{(n)}\right]$ as follows:

$$\sin^2\frac{\pi}{2\left(y_1^{(n)}-y_0^{(n)}\right)}\left(y-y_0^{(n)}\right),\quad \cos^2\frac{\pi}{2\left(y_1^{(n)}-y_0^{(n)}\right)}\left(y-y_0^{(n)}\right)$$

They have some interesting properties: one is that they can have the form of Adam probability essence wave functions; another is that they can also have the form of Eve probability essence wave functions; they are all at ground state, and they are all regarded as the probability essence wave functions of the descendant particles of $M$ when the quantum number is the natural number $n$.

When $x\in\left[x_1^{(n)},x_2^{(n)}\right]$, from the information of Adam and Eve probability essence wave functions $\sin^2(n\pi x)$ and $\cos^2(n\pi x)$ at the nodes as the following: $y_1^{(n)}=f\left(x_1^{(n)}\right), y_2^{(n)}=f\left(x_2^{(n)}\right)$, we can get Adam and Eve probability essence wave functions in $\left[y_1^{(n)},y_2^{(n)}\right]$ as follows:

$$\sin^2\frac{\pi}{2\left(y_2^{(n)}-y_1^{(n)}\right)}\left(y-y_1^{(n)}\right),\quad \cos^2\frac{\pi}{2\left(y_2^{(n)}-y_1^{(n)}\right)}\left(y-y_1^{(n)}\right)$$

They also have the properties: one is that they can have the form of Adam probability essence wave functions; another is that they can also have the form of Eve probability essence wave functions; they are all at ground state, and they are all regarded as the probability essence wave functions of the descendant particles of $M$ when the quantum number is also the natural number $n$.

At last, when $x\in\left[x_{m-1}^{(n)},x_m^{(n)}\right]$, from the information of Adam and Eve probability essence wave functions $\sin^2(n\pi x)$ and $\cos^2(n\pi x)$ at the nodes: $y_{m-1}^{(n)}=f\left(x_{m-1}^{(n)}\right), y_m^{(n)}=f\left(x_m^{(n)}\right)$, we get Adam and Eve probability essence wave functions in $\left[y_{m-1}^{(n)},y_m^{(n)}\right]$ as follows:

$$\sin^2\frac{\pi}{2\left(y_m^{(n)}-y_{m-1}^{(n)}\right)}\left(y-y_{m-1}^{(n)}\right),\quad \cos^2\frac{\pi}{2\left(y_m^{(n)}-y_{m-1}^{(n)}\right)}\left(y-y_{m-1}^{(n)}\right)$$

They have the same properties: one is that they can have the form of Adam probability essence wave functions; another is that they can also have the form of Eve probability essence wave functions; they are all at ground state, and they are all regarded as the probability essence wave functions of the descendant particles of $M$ when quantum number is $n$.

Based on these probability essence wave functions, we try to make some useful base functions defined respectively on the intervals $[0,1]$ and $[c,d]=\left[y_0^{(n)},y_m^{(n)}\right]$ denoted as follows:

$$A_0^{(n)}, A_1^{(n)}, \cdots, A_m^{(n)}, B_0^{(n)}, \quad B_1^{(n)}, \cdots, B_m^{(n)}:$$

$$A_0^{(n)}(x) = \chi_{\left[0,\frac{1}{m}\right]}(x)\cos^2(n\pi x),$$

$$A_1^{(n)}(x) = \chi_{\left[0,\frac{2}{m}\right]}(x)\sin^2(n\pi x),$$

$$A_2^{(n)}(x) = \chi_{\left[\frac{1}{m},\frac{3}{m}\right]}(x)\cos^2(n\pi x)$$

……

$$A_{m-2}^{(n)}(x) = \chi_{\left[\frac{m-3}{m},\frac{1}{m}\right]}(x)\cos^2(n\pi x),$$

$$A_{m-1}^{(n)}(x) = \chi_{\left[\frac{m-2}{m},1\right]}(x)\sin^2(n\pi x),$$

$$A_m^{(n)}(x) = \chi_{\left[\frac{m-1}{m},1\right]}(x)\cos^2(n\pi x);$$

$$B_0^{(n)}(y) = \chi_{\left[y_0^{(n)},y_1^{(n)}\right]}(y)\cos^2 \frac{\pi}{2\left(y_1^{(n)} - y_0^{(n)}\right)}\left(y - y_0^{(n)}\right),$$

$$B_1^{(n)}(y) = \chi_{\left[y_0^{(n)},y_1^{(n)}\right]}(y)\sin^2 \frac{\pi}{2\left(y_1^{(n)} - y_0^{(n)}\right)}\left(y - y_0^{(n)}\right)$$

$$+ \chi_{\left[y_1^{(n)},y_2^{(n)}\right]}(y)\cos^2 \frac{\pi}{2\left(y_2^{(n)} - y_1^{(n)}\right)}\left(y - y_1^{(n)}\right),$$

……

$$B_{m-1}^{(n)}(y) = \chi_{\left[y_{m-2}^{(n)},y_{m-1}^{(n)}\right]}(y)\sin^2 \frac{\pi}{2\left(y_{m-1}^{(n)} - y_{m-2}^{(n)}\right)}\left(y - y_{m-2}^{(n)}\right)$$

$$+ \chi_{\left[y_{m-1}^{(n)},y_m^{(n)}\right]}(y)\cos^2 \frac{\pi}{2\left(y_m^{(n)} - y_{m-1}^{(n)}\right)}\left(y - y_{m-1}^{(n)}\right),$$

$$B_m^{(n)}(y) = \chi_{\left[y_{m-1}^{(n)},y_m^{(n)}\right]}(y)\sin^2 \frac{\pi}{2\left(y_m^{(n)} - y_{m-1}^{(n)}\right)}\left(y - y_{m-1}^{(n)}\right)$$

where $\chi_A$ is the characteristic function of the set $A$; for example,

$$\chi_{\left[0,\frac{1}{m}\right]}(x) = \begin{cases} 1, & x \in \left[0, \dfrac{1}{m}\right], \\ 0, & x \in [0,1] - \left[0, \dfrac{1}{m}\right], \end{cases}$$

$$\chi_{\left[y_0^{(n)},y_1^{(n)}\right]}(y) = \begin{cases} 1, & y \in \left[y_0^{(n)}, y_1^{(n)}\right], \\ 0, & y \in [c,d] - \left[y_0^{(n)}, y_1^{(n)}\right] \end{cases}$$

Let us denote two classes of sets as the following:

$$\mathcal{A}(n) = \left\{A_0^{(n)}, A_1^{(n)}, \cdots, A_m^{(n)}\right\}, \quad \mathcal{B}(n) = \left\{B_0^{(n)}, B_1^{(n)}, \cdots, B_m^{(n)}\right\}$$

Clearly they are just the groups of base functions defined respectively on the closed intervals $X = [0,1]$ and $Y = [c,d]$. Clearly $\mathcal{A}(n)$ is a linearly independent group in the continuous function space $C[0,1]$ and $\mathcal{B}(n)$ is a linearly independent group in the continuous function space $C[c,d]$. Put

$$\mathcal{A}(n)\cdot\mathcal{B}(n)=\left\{A_i^{(n)}\cdot B_j^{(n)}\middle| i,j=0,1,\cdots,m\right\},$$
$$A_i^{(n)}\cdot B_j^{(n)}:[0,1]\times[c,d]\to[0,1]$$
$$(x,y)\mapsto\left(A_i^{(n)}\cdot B_j^{(n)}\right)(x,y)=A_i^{(n)}(x)\cdot B_j^{(n)}(y)$$

It is easy to know that $\mathcal{A}(n)\cdot\mathcal{B}(n)$ a linearly independent group in the continuous function space $C\left([0,1]\times[c,d]\right)$. Now we take the diagonal elements of $\mathcal{A}(n)\cdot\mathcal{B}(n)$ to make a set as follows:

$$\mathcal{C}(n)=\left\{A_i^{(n)}\cdot B_i^{(n)}\middle| i=0,1,\cdots,m\right\},$$

which is clearly a linearly independent group with $m+1=2n+1$ dimension in the continuous function space $C\left([0,1]\times[c,d]\right)$. By using $\mathcal{C}(n)$, we can get a sequence of binary nonnegative continuous functions as the following:

$$\mu_n:[0,1]\times[c,d]\to[0,1]$$
$$(x,y)\mapsto\mu_n(x,y)=\mathop{\vee}_{i=1}^{m}\left[A_i^{(n)}(x)\cdot B_i^{(n)}(y)\right],\quad n=1,2,3,\cdots$$

where

$$\mathop{\vee}_{i=0}^{m}\left[A_i^{(n)}(x)\cdot B_i^{(n)}(y)\right]=\max_{0\le i\le m}\left\{A_i^{(n)}(x)\cdot B_i^{(n)}(y)\right\}.$$

Then this sequence of binary nonnegative continuous functions as being $\{\mu_n(x,y)\}_{n=1}^{\infty}$ are normalized as the following:

$$p_n(x,y)=\chi_{[0,1]\times[c,d]}(x,y)\frac{\mu_n(x,y)}{\int_c^d\int_0^1\mu_n(x,y)dxdy},\quad n=1,2,3,\cdots,$$

where

$$\chi_{[0,1]\times[c,d]}(x,y)=\begin{cases}1, & (x,y)\in[0,1]\times[c,d],\\ 0, & (x,y)\in\mathbb{R}^2-[0,1]\times[c,d]\end{cases}$$

Therefore $\{p_n(x,y)\}_{n=1}^{\infty}$ becomes a sequence of probability density functions defined on $\mathbb{R}^2$, and $p_n(x,y)$ is called the probability density function when the quantum number is just $n$.

And now by means of the sequence $\{p_n(x,y)\}_{n=1}^{\infty}$, we can construct a sequence of functions of one variable as follows:

$$f_n(x)=\frac{\int_{-\infty}^{+\infty}yp_n(x,y)dy}{\int_{-\infty}^{+\infty}p_n(x,y)dy},\quad n=1,2,3,\cdots \tag{4.18}$$

Apparently, $\{f_n(x)\}_{n=1}^{\infty}$ is just the sequence of conditional mathematical expectations formed by $\{p_n(x,y)\}_{n=1}^{\infty}$.

**Case 2**. Suppose $f:X\to Y$ be not strict monotonous function and not constant function.

Because the elements of the set $Y(n)$ may not always satisfy the monotonicity about the subscript $i$ as like as $y_0\le y_1\le\cdots\le y_m$, it is of a little difficulty to make the continuous base functions as follows:

$$B_i^{(n)}(y),\quad i=0,1,\cdots,m.$$

So we have to make a permutation on the subscript set $\{0,1,\cdots,m\}$ as the following:

$$\sigma=\begin{pmatrix}0 & 1 & \cdots & m\\ k_0 & k_1 & \cdots & k_m\end{pmatrix},$$
$$\left(\forall i\in\{0,1,\cdots,m\}\right)\left(k_i=\sigma(i)\right)$$

such that the subscript set after the permutation is denoted by the following symbol:

$$K(n)=\{k_0,k_1,\cdots,k_m\}$$

and satisfies the following condition:

$$c(n)=y_{k_0}^{(n)}\le y_{k_1}^{(n)}\le\cdots\le y_{k_m}^{(n)}=d(n). \tag{4.19}$$

Since (4.19) shows that the inequalities may not be strict, we have to consider the following two situations.

1) Assume that $c(n)=y_{k_0}^{(n)}<y_{k_1}^{(n)}<\cdots<y_{k_m}^{(n)}=d(n)$. Based on these nodes $y_{k_0}^{(n)},y_{k_1}^{(n)},\cdots,y_{k_m}^{(n)}$ in $[c,d]$ and doing in imitation of Case 1, we can get the continuous base functions $B_{k_j}^{(n)}$ as the following:

$$B_{k_0}^{(n)}(y)=\chi_{\left[y_{k_0}^{(n)},y_{k_1}^{(n)}\right]}(y)\cos^2\frac{\pi}{2\left(y_{k_1}^{(n)}-y_{k_0}^{(n)}\right)}\left(y-y_{k_0}^{(n)}\right),$$

$$B_{k_1}^{(n)}(y)=\chi_{\left[y_{k_0}^{(n)},y_{k_1}^{(n)}\right]}(y)\sin^2\frac{\pi}{2\left(y_{k_1}^{(n)}-y_{k_0}^{(n)}\right)}\left(y-y_{k_0}^{(n)}\right)$$
$$+\chi_{\left[y_{k_1}^{(n)},y_{k_2}^{(n)}\right]}(y)\cos^2\frac{\pi}{2\left(y_{k_2}^{(n)}-y_{k_1}^{(n)}\right)}\left(y-y_{k_1}^{(n)}\right),$$

……

$$B_{k_{m-1}}^{(n)}(y)=\chi_{\left[y_{k_{m-2}}^{(n)},y_{k_{m-1}}^{(n)}\right]}(y)\sin^2\frac{\pi}{2\left(y_{k_{m-1}}^{(n)}-y_{k_{m-2}}^{(n)}\right)}\left(y-y_{k_{m-2}}^{(n)}\right)$$
$$+\chi_{\left[y_{k_{m-1}}^{(n)},y_{k_m}^{(n)}\right]}(y)\cos^2\frac{\pi}{2\left(y_{k_m}^{(n)}-y_{k_{m-1}}^{(n)}\right)}\left(y-y_{k_{m-1}}^{(n)}\right),$$

$$B_{k_m}^{(n)}(y)=\chi_{\left[y_{k_{m-1}}^{(n)},y_{k_m}^{(n)}\right]}(y)\sin^2\frac{\pi}{2\left(y_{k_m}^{(n)}-y_{k_{m-1}}^{(n)}\right)}\left(y-y_{k_{m-1}}^{(n)}\right)$$

Then we easily make a sequence of binary nonnegative continuous functions defined on $X\times Y=[0,1]\times[c,d]$ as follows:

$$\mu_n(x,y)=\bigvee_{j=1}^{m}\left[A_{k_j}^{(n)}(x)\cdot B_{k_j}^{(n)}(y)\right],\quad n=1,2,3,\cdots \tag{4.20}$$

2) Assume $c(n)=y_{k_0}^{(n)}\le y_{k_1}^{(n)}\le\cdots\le y_{k_m}^{(n)}=d(n)$. Firstly we do a kind of screen work on the elements in the following node set:

$$Y(n)=\left\{y_{k_0}^{(n)},y_{k_1}^{(n)},\cdots,y_{k_m}^{(n)}\right\}.$$

In fact, let $K(n)=\{k_0,k_1,\cdots,k_m\}$. We define an equivalence relation on the set $K(n)$ as being "$\sim$" as follows:

$$\left(\forall s,t\in\{0,1,\cdots,m\}\right)\left(k_s\sim k_t\Leftrightarrow y_{k_s}^{(n)}=y_{k_t}^{(n)}\right).$$

Then we get the quotient set of $K(n)$ as the following:

$$K(n)/\sim=\left\{\left[k_j\right]\middle|\, j=0,1,\cdots,m\right\},$$

where $\left[k_j\right]$ is the equivalence class in which $k_j$ belongs.

Let all the elements of the quotient set $K(n)/\sim$ be the following:

$$\left[k_{j_0}\right],\left[k_{j_1}\right],\cdots,\left[k_{j_{q(m)}}\right],$$

where $0\le q(m)\le m$, and stipulate the representative element $k_{j_s}$ be the smallest element in $\left[k_{j_s}\right]$. Thus we have the following inequalities:

$$y_{k_{j_0}}^{(n)} < y_{k_{j_1}}^{(n)} < \cdots < y_{k_{j_{q(m)}}}^{(n)} .$$

Based on the nodes $y_{k_{j_0}}^{(n)}, y_{k_{j_1}}^{(n)}, \cdots, y_{k_{j_{q(m)}}}^{(n)}$ in $[c,d]$, we make the continuous base functions $B_{k_{j_s}}^{(n)} \left(s = 0,1,\cdots,q(m)\right)$ as follows:

$$B_{k_{j_0}}^{(n)}(y) = \chi_{\left[y_{k_{j_0}}^{(n)}, y_{k_{j_1}}^{(n)}\right]}(y) \cos^2 \frac{\pi}{2\left(y_{k_{j_1}}^{(n)} - y_{k_{j_0}}^{(n)}\right)}\left(y - y_{k_{j_0}}^{(n)}\right),$$

$$B_{k_{j_1}}^{(n)}(y) = \chi_{\left[y_{k_{j_0}}^{(n)}, y_{k_{j_1}}^{(n)}\right]}(y) \sin^2 \frac{\pi}{2\left(y_{k_1}^{(n)} - y_{k_0}^{(n)}\right)}\left(y - y_{k_{j_0}}^{(n)}\right)$$

$$+ \chi_{\left[y_{k_{j_1}}^{(n)}, y_{k_{j_2}}^{(n)}\right]}(y) \cos^2 \frac{\pi}{2\left(y_{k_{j_2}}^{(n)} - y_{k_{j_1}}^{(n)}\right)}\left(y - y_{k_{j_1}}^{(n)}\right),$$

……

$$B_{k_{j_{q(m)-1}}}^{(n)}(y) = \chi_{\left[y_{k_{j_{q(m)-2}}}^{(n)}, y_{k_{j_{q(m)-1}}}^{(n)}\right]}(y) \sin^2 \frac{\pi}{2\left(y_{k_{j_{q(m)-1}}}^{(n)} - y_{k_{j_{q(m)-2}}}^{(n)}\right)}\left(y - y_{k_{j_{q(m)-2}}}^{(n)}\right)$$

$$+ \chi_{\left[y_{k_{j_{q(m)-1}}}^{(n)}, y_{k_{j_{q(m)}}}^{(n)}\right]}(y) \cos^2 \frac{\pi}{2\left(y_{k_{j_{q(m)}}}^{(n)} - y_{k_{j_{q(m)-1}}}^{(n)}\right)}\left(y - y_{k_{j_{q(m)-1}}}^{(n)}\right),$$

$$B_{k_{j_{q(m)}}}^{(n)}(y) = \chi_{\left[y_{k_{j_{q(m)-1}}}^{(n)}, y_{k_{j_{q(m)}}}^{(n)}\right]}(y) \sin^2 \frac{\pi}{2\left(y_{k_{j_{q(m)}}}^{(n)} - y_{k_{j_{q(m)-1}}}^{(n)}\right)}\left(y - y_{k_{j_{q(m)-1}}}^{(n)}\right)$$

Hence for the nodes $y_{k_{j_0}}^{(n)}, y_{k_{j_1}}^{(n)}, \cdots, y_{k_{j_{q(m)}}}^{(n)}$ which correspond to the representative elements:

$$k_{j_0}, k_{j_1}, \cdots, k_{j_{q(m)}}$$

coming from these equivalence classes $\left[k_{j_0}\right], \left[k_{j_1}\right], \cdots, \left[k_{j_{q(m)}}\right]$, we have made the continuous base functions as follows:

$$B_{k_{j_0}}^{(n)}(y), B_{k_{j_1}}^{(n)}(y), \cdots, B_{k_{j_{q(m)}}}^{(n)}(y) .$$

For any $s \in \left\{0,1,\cdots,q(m)\right\}$ and we can define the continuous base functions corresponding to the elements in $\left[k_{j_s}\right] - \left\{k_{j_s}\right\}$ as following:

$$\left(\forall \tau \in \left[k_{j_s}\right] - \left\{k_{j_s}\right\}\right)\left(B_{\tau}^{(n)}(y) \equiv B_{k_{j_s}}^{(n)}(y)\right)$$

So for all the nodes $y_{k_0}^{(n)} \le y_{k_1}^{(n)} \le \cdots \le y_{k_m}^{(n)}$ in $[c,d]$, we have got the corresponding continuous base functions as follows:

$$B_{k_0}^{(n)}(y), B_{k_1}^{(n)}(y), \cdots, B_{k_m}^{(n)}(y) .$$

By using these continuous base functions, we get a sequence of binary nonnegative continuous functions defined on $X \times Y = [0,1] \times [c,d]$ as the following:

$$\mu_n(x,y) = \mathop{\vee}_{j=1}^{m} \left[ A_{k_j}^{(n)}(x) \cdot B_{k_j}^{(n)}(y) \right], \quad n = 1,2,3,\cdots \tag{4.21}$$

Based on above two cases, we have got the sequence of binary nonnegative continuous functions defined on $X \times Y = [0,1] \times [c,d]$ as being $\left\{\mu_n(x,y)\right\}_{n=1}^{\infty}$. Now we normalize $\left\{\mu_n(x,y)\right\}_{n=1}^{\infty}$ as follows:

$$p_n(x,y)=\chi_{[0,1]\times[c,d]}(x,y)\frac{\mu_n(x,y)}{\int_c^d\int_0^1\mu_n(x,y)dxdy},$$
$$n=1,2,3,\cdots$$

where

$$\chi_{[0,1]\times[c,d]}(x,y)=\begin{cases}1, & (x,y)\in[0,1]\times[c,d],\\ 0, & (x,y)\in\mathbb{R}^2-[0,1]\times[c,d]\end{cases}$$

Therefore $\{p_n(x,y)\}_{n=1}^{\infty}$ becomes a sequence of probability density functions defined on $\mathbb{R}^2$, and $p_n(x,y)$ is also called the probability density function when the quantum number is just $n$. And by means of $\{p_n(x,y)\}_{n=1}^{\infty}$, we can construct a sequence of functions of one variable defined on $[0,1]$ as follows:

$$f_n(x)=\frac{\int_{-\infty}^{+\infty}yp_n(x,y)dy}{\int_{-\infty}^{+\infty}p_n(x,y)dy},\quad x\in[0,1], \tag{4.22}$$
$$n=1,2,3,\cdots$$

Apparently, $\{f_n(x)\}_{n=1}^{\infty}$ is just the sequence of conditional mathematical expectations formed by $\{p_n(x,y)\}_{n=1}^{\infty}$. Besides, it is not under the following expression:

$$f_n(x)=\frac{\int_{-\infty}^{+\infty}y\mu_n(x,y)dy}{\int_{-\infty}^{+\infty}\mu_n(x,y)dy},$$
$$x\in[0,1],\quad n=1,2,3,\cdots$$

**Step 7**. Similar to the proof of theorem 2.1, we know that the sequence of conditional mathematical expectations $\{f_n(x)\}_{n=1}^{\infty}$ can uniformly converge to $f(x)$ on the closed interval $[0,1]$.

Paying attention to the process of the theorem, when the quantum number is $n$, the set of the descendant particles generated by the particle $M$ is the following:

$$\mathcal{M}_n=\{M_1^{(n)},M_2^{(n)},\cdots,M_{2n}^{(n)}\},$$

When $n\to\infty$, the set of all descendant particles generated by the particle $M$ is $\mathcal{M}=\bigcup_{n=1}^{\infty}\mathcal{M}_n$. Clearly the cardinal number of the set is as being: $\text{card}(\mathcal{M})=\aleph_0$; i.e., we all use countable infinite particles. These particles can be expressed as the following expression:

$$M\overset{n=1,2,3,\cdots}{\Rightarrow}\begin{cases}M_1^{(1)},M_2^{(1)};\\ M_1^{(2)},M_2^{(2)},M_3^{(2)},M_4^{(2)};\\ M_1^{(3)},M_2^{(3)},M_3^{(3)},M_4^{(3)},M_5^{(3)},M_6^{(3)};\\ \cdots\cdots\\ M_1^{(n)},M_2^{(n)},\cdots,M_{2n-1}^{(n)},M_{2n}^{(n)};\\ \cdots\cdots\end{cases}$$

where only the particle $M$ moves along $Ox$ axis, but all the descendant particles $M_1^{(1)},M_2^{(1)},\cdots M_1^{(n)},\cdots M_{2n}^{(n)},\cdots$ move along $Oy$ axis.

This means that the motion curve of a mass point in classic physics $y=f(x)$ can be constructed by an infinite sequence of microscopic particles wave functions. In other words, this motion curve of a mass point $y=f(x)$ have been quantization, which is the limit state of these microscopic particles wave functions when $n \to \infty$. Clearly this fact meets the Bohr's correspondence principle.
We finally end the proof of the theorem. □

**Example 4.1** Suppose we cast an object $B$ with quality $m_0$, which is regarded as a mass point. So the movement of $B$ can be described by its equation of locus as follows:

$$y = f(x) = x\tan\alpha - \frac{g}{2v_0^2\cos^2\alpha}x^2,$$

$$x \in [0, d_0], \quad d_0 = \frac{v_0^2}{g}\sin 2\alpha$$

where $\alpha \in \left(0, \frac{\pi}{2}\right)$ is a mass ejection angle, $d_0 \in (0, +\infty)$ is the maximum range of fire, and $v_0 \in (0, +\infty)$ is the initial velocity; here the air friction is omitted. Clearly $y(x) \in C[0, d_0]$, which means that the projectile motion is expressed by a unary continuous function.

Now if we take $\alpha = \frac{\pi}{4}, v_0 = \sqrt{g}$, then $d_0 = 1$; then we have the following equation:

$$y = f(x) = x - x^2 = x(1-x).$$

When the quantum number $n = 5, 10, 20$, the approximation situations of the sequence of conditional mathematical expectations $f_n(x)$ to $f(x)$ are respectively shown in Figs, 4.5, 4.6 and 4.7., where red curve means $f_n(x)$, and blue curve indicates $f(x)$. □

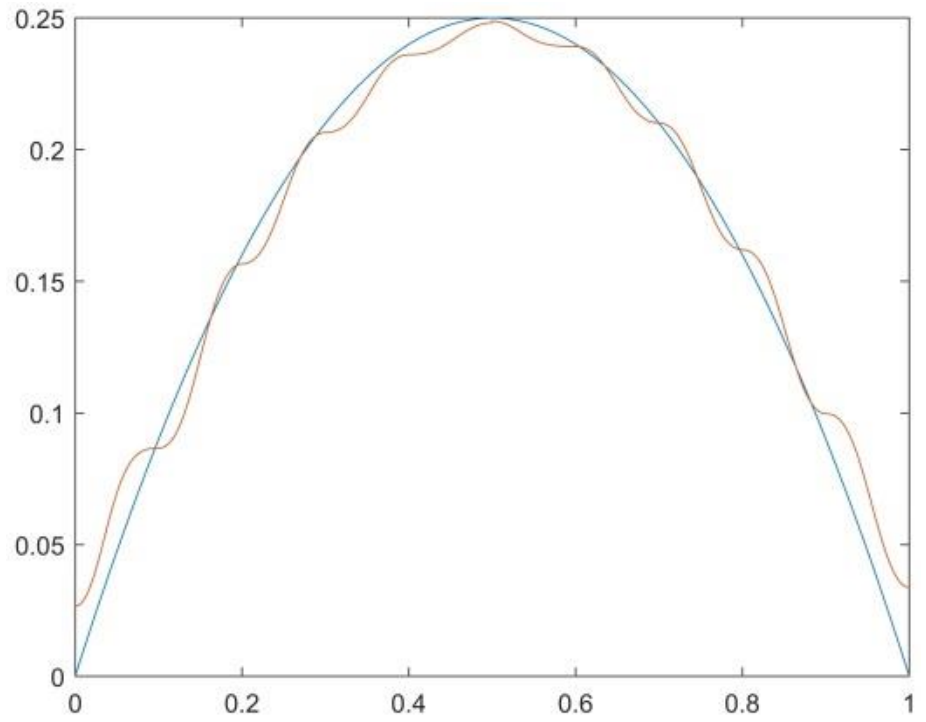

Fig. 4.5. Approximation of $f_5(x)$ to $f(x)$

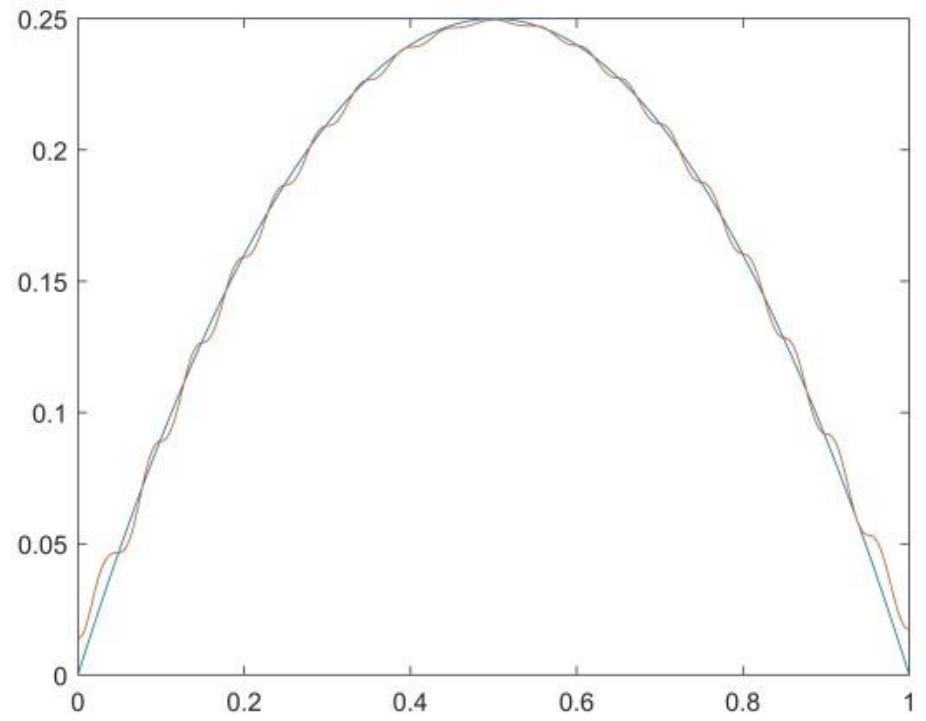

Fig. 4.6. Approximation of $f_{10}(x)$ to $f(x)$

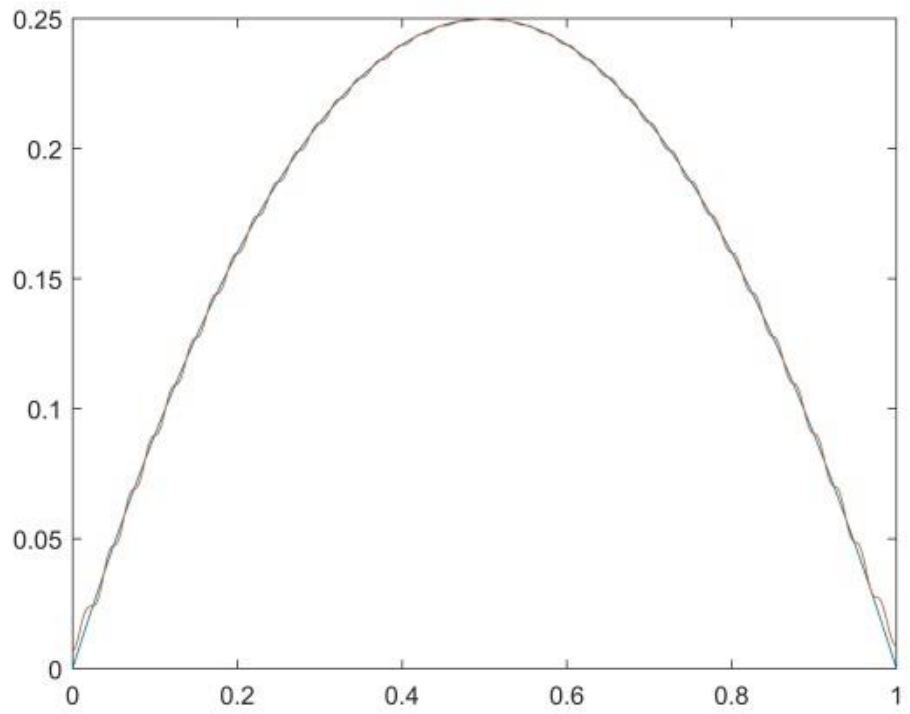

Fig. 4.7. Approximation of $f_{20}(x)$ to $f(x)$

## 5. Duality of Mass Point Motion

We firstly review the projectile motion in Example 4.1. The property of mass point motion is shown as its momentum $p = m_0 v_0$ and its energy as the following:

$$E = E_k = \frac{1}{2} m_0 v_0^2 .$$

Actually, more straightway, its property of mass point should be described by its equation of locus as the following:

$$y = f(x) = x \tan\alpha - \frac{g}{2v_0^2 \cos^2\alpha} x^2$$

In other words, the property of mass point can be described by its momentum and energy or by its equation of locus; these two methods are equivalent.

Then we ask an interesting and important problem: is there wave nature on mass point motion in classic physics? Alternatively, we can ask the question: is there wave mass point duality in classic physics?

For answering this problem, we firstly review the particle nature and wave nature in quantum mechanics. As we all know, a microscopic particle has no determinate movement locus so that it has no an equation describing its movement locus. Thus, its nature of particle can only be described by its momentum $p = mv$ and its energy $E = \frac{1}{2} mv^2$. Based on the viewpoint of de Broglie, an object particle is of wave-particle duality, which means the particle also has its nature of wave. The nature of wave should be shown by its wave function $\Psi$, and the wave function $\Psi$ should be the solution of Schrodinger Equation. The wave as being the solution of Schrodinger Equation is called de Broglie wave. Then Born gave Schrodinger Equation the statistical interpretation of de Broglie wave, which means that $|\Psi|^2$ should be a kind of probability density function. So $|\Psi|^2$ is often called probability wave. In fact, in quantum mechanics, the probability wave $|\Psi|^2$ is much more important than the wave function $\Psi$ itself.

Again, we consider the movement of the particle in the infinite deep square potential well as we have discussed in Step 1 in Theorem 2.1, where the wave function is as following:

$$\psi_n(x) = \sqrt{2} \sin(n\pi x),$$
$$x \in [0,1], \quad n = 1, 2, 3, \cdots$$

Then, its probability wave is $|\psi_n(x)|^2 = 2\sin^2(n\pi x)$, which figure is shown in Fig. 5.1.

It is worth noting that, the probability wave $|\psi_n(x)|^2$ describes the probability density that the particle $M$ appears at $x$ in $[0,1]$ when the quantum number is $n$. Because the particle $M$ does one-dimension motion along $Ox$ axis, $|\psi_n(x)|^2$ is a curve on two-dimension plane.

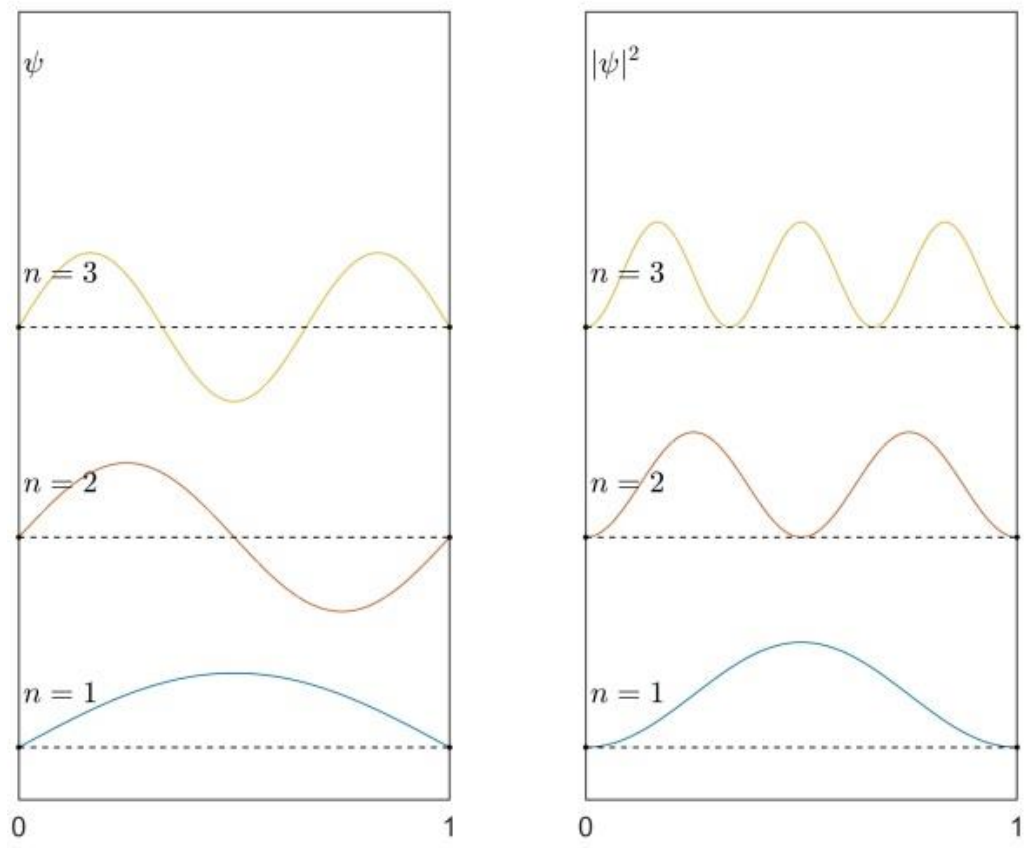

Fig. 5.1. The nature of waves of $\psi_n(x)$ and $|\psi_n(x)|^2$

It is well-known that the wave nature of simple harmonic wave is constructed by its frequency $\nu$ and its wave length $\lambda$. When the quantum number is $n$, its energy expression is $E_n = \frac{n^2\pi^2\hbar^2}{2m}$, and the wave frequency is as following:

$$\nu_n = \frac{1}{T_n} = \frac{k_n}{2\pi} = \frac{n\pi}{2\pi} = \frac{n}{2} = \frac{\sqrt{mE_n}}{\sqrt{2}\pi\hbar}.$$

Based on the definition of wave length, we know the wave length is $\lambda_n = \frac{2}{n}$ so that

$$\lambda_n = \frac{2}{n} = \frac{\sqrt{2}\pi\hbar}{\sqrt{mE_n}}.$$

This just gives the result that $\nu_n \cdot \lambda_n = 1$, which means that the relation between the wave nature and the particle nature can be established by using Planck number $\hbar$.

Now we return to continue to discuss the motion of projectile. Its mass point nature reflected in its equation of locus.

Especially, when $\alpha = \frac{\pi}{4}, v_0 = \sqrt{g}$, the equation of locus is as follows:

$$y = f(x) = x - x^2 = x(1-x), \quad x \in [0,1].$$

Because this sequence of conditional mathematical expectations as being $\{f_n(x)\}_{n=1}^{\infty}$ uniformly converges to $y = f(x)$ in $[0,1]$, for arbitrarily given a $\varepsilon > 0$, there must exist a natural number $N \in \mathbb{N}_+$, such that

$$(\forall n \in \mathbb{N}_+)(n > N \Rightarrow \|f_n - f\| < \varepsilon),$$

where $\|\cdot\|$ is a kind of norm in the linear normed space $(C[0,1], \|\cdot\|)$ and defined as the following:

$$(\forall f \in C[0,1])\left(\|f\| = \max_{x\in[0,1]} |f(x)|\right).$$

For $\varepsilon > 0$ is small enough, that $\|f_n - f\| < \varepsilon$ means that the difference between $f_n$ and $f$ is very small so that $f_n$ can be replaced by $f$ approximately.

We now take notice of the following important expression:

$$f(x) \approx f_n(x) = \frac{\int_c^d y p_n(x,y)dy}{\int_c^d p_n(x,y)dy} = \frac{\int_0^1 y p_n(x,y)dy}{\int_0^1 p_n(x,y)dy},$$

for above the motion of projectile where $c=0, d=1$, where $p_n(x,y)$ is a binary probability density function.

When the quantum number $n=5,10,15$, the graphs of the probability density function $p_n(x,y)$ are respectively shown in Figs. 5.2, 5.3 and 5.4.

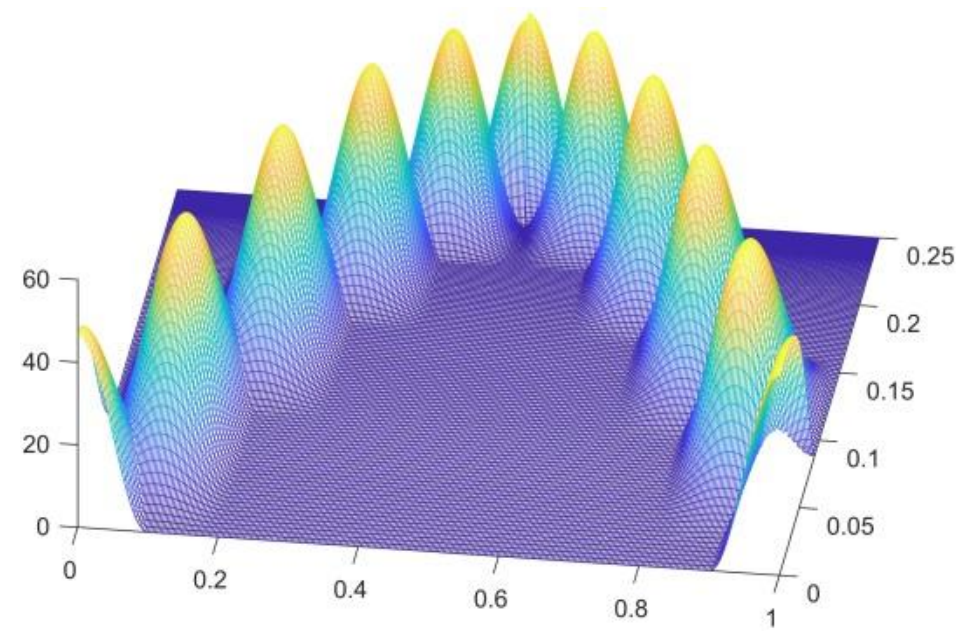

Fig. 5.2. Graph of $p_5(x,y)$

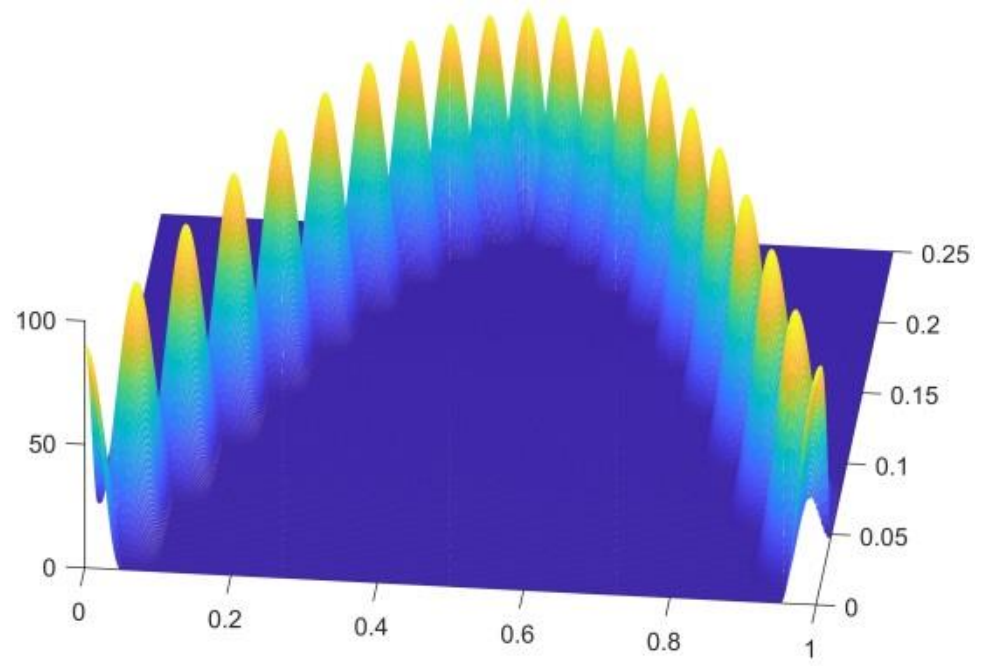

Fig. 5.3. Graph of $p_{10}(x,y)$

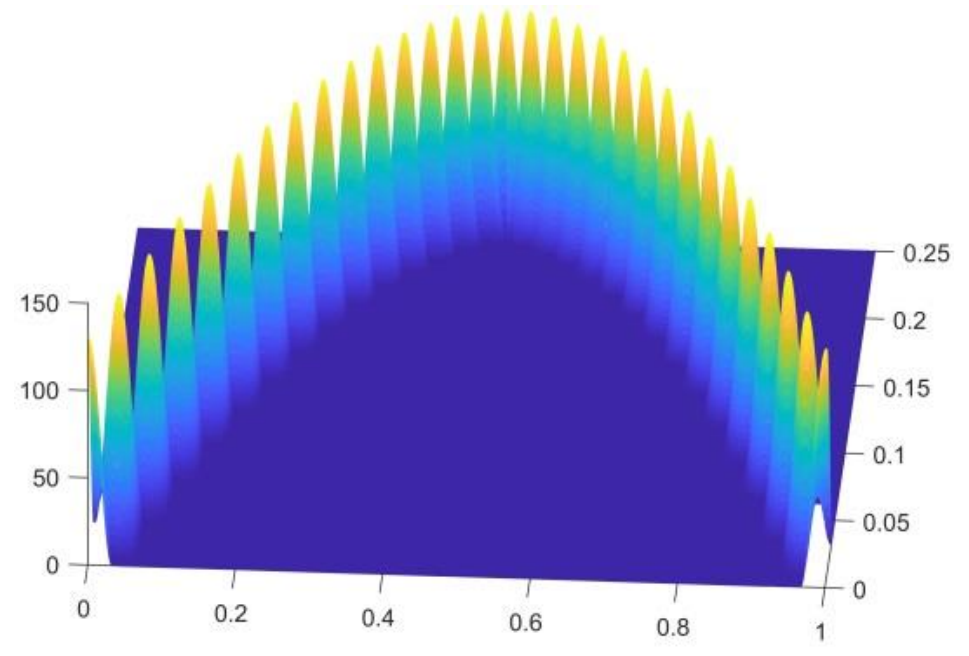

Fig. 5.4. Graph of $p_{15}(x,y)$

Apparently, the probability density function $p_n(x,y)$ shows up waviness. We observe the motion curve of the projectile, and suppose some mass point $B$ moves in the rectangle as in Fig 5.5, and the probability density function that $B$ falls into the set of graph of $f(x)$ as follows

$$G_f = \{(x,y) \in [0,1]\times[0,0.25] \mid y = f(x)\} \tag{5.1}$$

is just $p_n(x,y)$.

It is worth noting that, since the mass point $B$ moves in a two-dimension region, the probability density function $p_n(x,y)$ is a wave surface in three-dimension space. From Fig. 5.1, we can learn that, since the particle $M$ moves in $[0,1]$ on $Ox$ axis, the probability wave $|\psi_n(x)|^2$ mainly roots in $[0,1]$; while from Fig. 3.4, we also can learn that, since the mass point $B$ moves in $G_f$ on $x-y$ plane, the probability wave $p_n(x,y)$ roots in $G_f$.

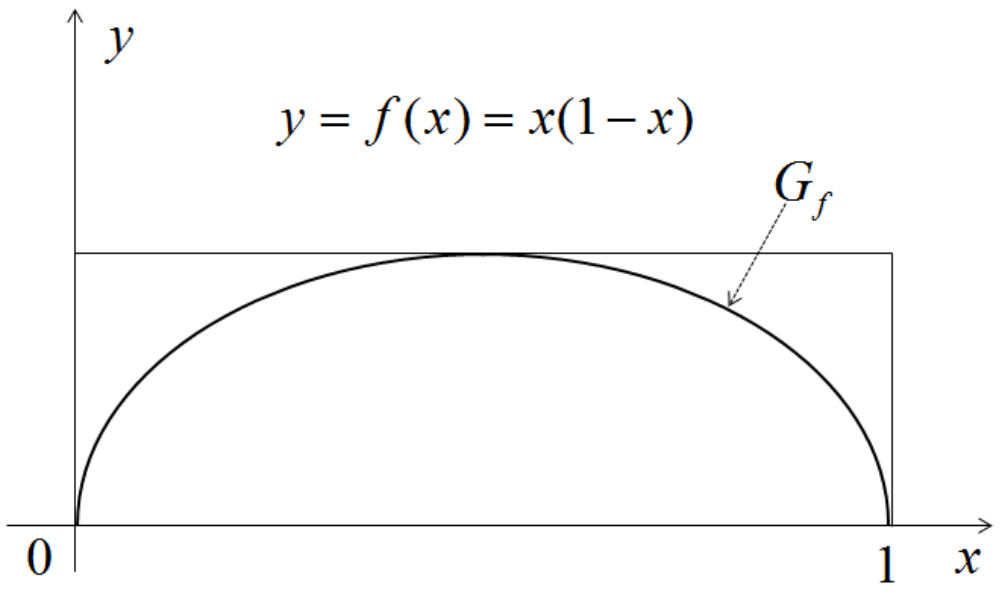

Fig. 5.5. Graph of the motion of projectile

Above discussion reveals an important conclusion: the motion of mass point in classic mechanics is surely of waviness so that the motion of mass point in classic mechanics also has wave mass point duality, which is very same with wave-particle duality in quantum mechanics.

Furthermore, the relationship between the wave nature and particle nature is established by means of Schrodinger Equation and the energy of the particle $E$ and the momentum of the particle $p$ can be respectively expressed by the frequency $\nu$ and the wavelength $\lambda$ of the particle as the following:

$$E = 2\pi\hbar\nu, \quad p = \frac{2\pi\hbar}{\lambda}.$$

While in classic mechanics, the relation between the mass point nature and waviness of motion of mass point is related by means of the following integral equation:

$$\frac{\int_c^d yp(x,y)dy}{\int_c^d p(x,y)dy} = f(x), \qquad (x,y) \in [a,b]\times[c,d] \tag{5.2}$$

where $p(x,y) \in C([a,b]\times[c,d])$ is an unknown binary function satisfying the following conditions:

(1) $(\forall(x,y)\in[a,b]\times[c,d])(p(x,y)\geq 0)$;

(2) $(\forall x\in[a,b])\left(\int_c^d p(x,y)\mathrm{d}y > 0\right)$.

(3) $\int_c^d \int_a^b p(x,y)\mathrm{d}x\mathrm{d}y = 1$.

Because $y = f(x)$ is the equation of locus of motion of the mass point, it completely represents the mass point nature of motion of the mass point; while $p(x,y)$ is the probability density function which is the probability wave of itself so that $p(x,y)$ itself represents the waviness of motion of the mass point. And the relation between the mass point nature and the wave nature is related by means of the integral equation (8.3.2). This adequately explains that the motion of mass point in classic mechanics has the duality of wave mass point, or written by **wave-mass-point duality**.

Here we need to explain that to solve the integral equation (5.2) is not an easy thing; however, we have given a kind of approximate method to do it; actually, $\{p_n(x,y)\}_{n=1}^{\infty}$ is a sequence of approximate solutions of the integral equation because if we write

$$f_n(x) = \frac{\int_c^d yp_n(x,y)dy}{\int_c^d p_n(x,y)dy},$$

then we have $\lim_{n\to\infty}\| f_n - f\| = 0$ based on Theorem 2.1.

## 6. Conclusions

Firstly from physical world, we can receive a point of view: continuous functions can describe a large proportion of certainty phenomena; for example, the trajectory of a projectile motion is just described as a continuous function. And random variables or vectors should describe random phenomena. So if we want to consider the connection between some certainty phenomenoa and some random phenomena, we should or must research the relation between continuous functions and random vectors. Theorem 2.1 shows us an interesting conclusion:

For arbitrarily given a continuous function $f(x)\in C[a,b]$, there must be a sequenc of probability spaces $\{(\Omega,\mathcal{F},P_n)\}$ and a sequence of random vectors $\{(\xi_n,\eta_n)\}$, where every random vector $(\xi_n,\eta_n)$ $(\xi_n,\eta_n)$ is defined on the probability space $(\Omega,\mathcal{F},P_n)$, such that the sequence of conditional mathematical expectations $\{E(\eta_n|\xi_n=x)\}$ converges uniformly to the continuous function $f(x)$ in $[a,b]$.

This is random vector representation of continuous functions, which is like a bridge to be set up between real function theory and probability theory. By using this conclusion, we have a result with respect to function approximation which has been shown by theorem 3.1:

For any continuous function $f(x)\in C[a,b]$, but assuming $f(x)$ not being constant function, if $\{E(\eta_n|\xi_n=x)\}$ is called a sequence of conditional methematical expectations generated by the continuous function $f(x)$, then by means of $\{E(\eta_n|\xi_n=x)\}$ we can construct a group of continuous functions:

$$\Phi(n)=\{\varphi_0^{(n)}(x),\varphi_1^{(n)}(x),\cdots,\varphi_n^{(n)}(x)\},$$

where $\varphi_l^{(n)}(x)\in C[a,b]$, such that, by using the sequence of the groups of base functions $\{\Phi(n)\}$, the sequence of interpolation functions constructed as the following

$$f_n(x)=\sum_{l=0}^{n}\varphi_l^{(n)}(x)y_l^{(n)},$$

$$n=1,2,3,\cdots$$

uniformly converges to $f(x)$ in $[a,b]$.

And then, in approximation from a sequence of random vectors to a continuous function, the base functions are appropriately selected by us, an important conclusion for quantum mechanics is deduced: classical mechanics and quantum mechanics is unified. This is the content of theorem 4.1:

Given arbitrarily a non-constant function $f(x)\in C[0,1]$, there must exist some microscopic particles such that the limit of the group behavior of these microscopic particles is just this continuous function $f(x)$ when the quantum number $n\to\infty$.

Particularly, an interesting and very important conclusion is introduced as the fact that the mass point motion of a macroscopical object possesses a kind of wave characteristic curve, which is called wave-mass-point duality.

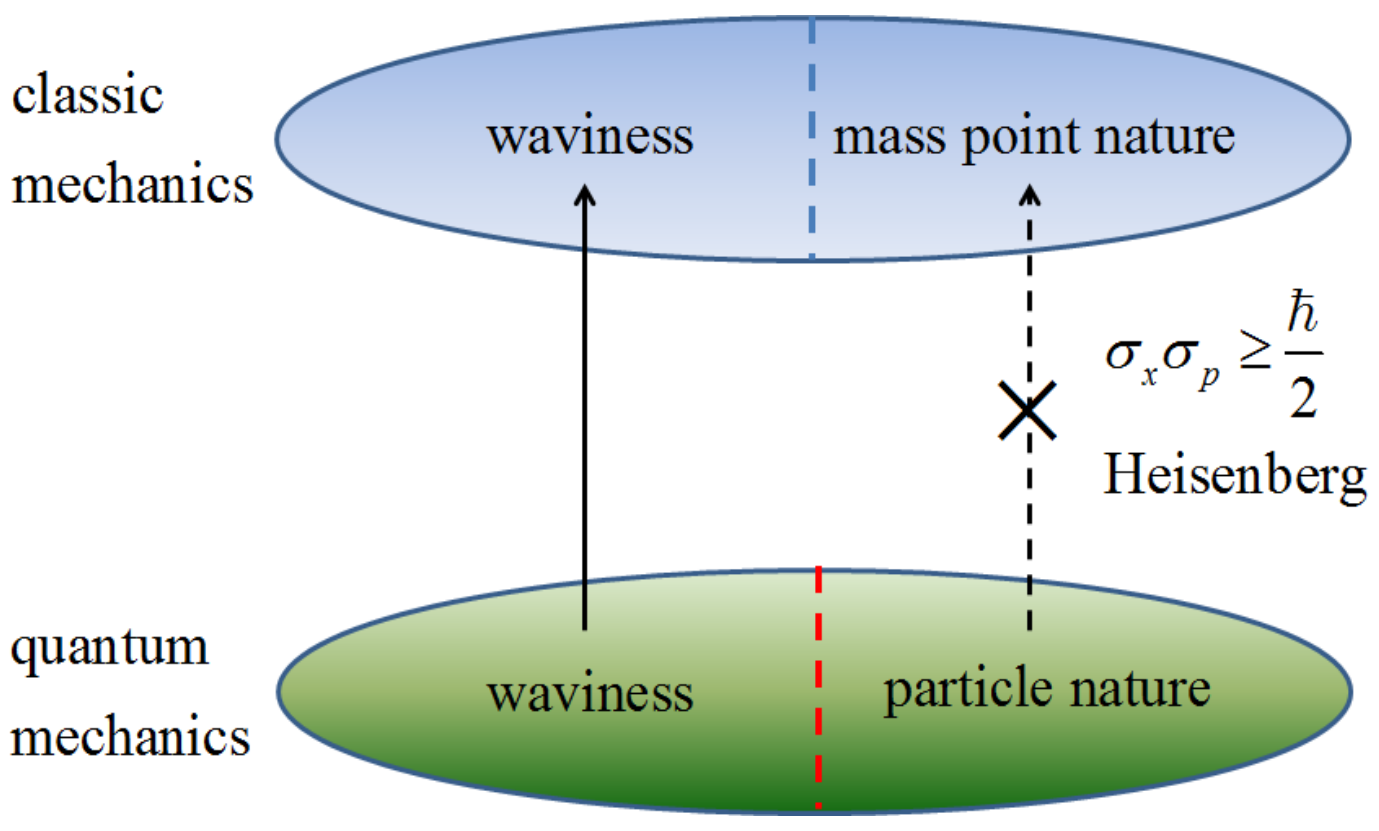

Fig. 6.1. Unified frame of two kinds of mechanics

In this paper, we reveal an important problem: unified theory of classic mechanics and quantum mechanics. So-called unified theory here means almost every motion of a mass point in classic mechanics can be represented by the motions of an infinite sequence of particles in quantum mechanics, where limit operation plays an important role in the unified theory. Clearly this situation is just according with Bohr's Correspondence Principle.

It is worth noting that this kind of correspondence relation between classic mechanics and quantum mechanics cannot be expressed by the relationship between the mass point nature in classic mechanics and the particle nature in quantum mechanics because of Heisenberg's Uncertainty Principle (see Fig. 6.1). As we all know, in classic mechanics, the motion of a mass point has no uncertainty so that we can use continuous functions to describe the movement locus of the mass point. However, in quantum mechanics, the motion of a particle has surely uncertainty so that we cannot use continuous functions to describe the movement locus of the particle. By now, we have known that the position and momentum of a particle are all random and they are related by Planck constant $\hbar$, i.e.,

$$\sigma_x \sigma_p \geq \frac{\hbar}{2}.$$

Fortunately, we have pointed that the motion of a mass point in classic mechanics has also waviness in Section 5. The wave function of the motion of a mass point has surely no uncertainty. On the other hand, although the motion of a particle has surely uncertainty, the wave function of the particle must have no uncertainty. Thus, we can consider the relation between the wave function of a mass point in classic mechanics and the wave functions of some particles in quantum mechanics. As we discussed in Section 4, we have revealed the relation by means of Theorem 4.1. In other words, by using wave functions of both classic mechanics and quantum mechanics, classic mechanics and quantum mechanics are unified, which is the significance of our unified theory about the two kinds of mechanics.

We need to emphasize our new and important and interesting conclusion: The motion of a mass point has also so-called duality: wave-mass-point duality, which is very similar to the case of the motion of a particle in quantum mechanics and is an important support to our unified theory on classic mechanics and quantum mechanics. It is not difficult to understand that Theorem 4.1 should be the most important in physics.

Prigogine had ever pointed out his conclusion by many experiments: world is random not certain (see [5]). In fact, Theorem 2.1 just prove his idea, because, as we all know, a large part of physical phenomenon can be described by some kind of continuous functions, and based on Theorem 2.1, any one of these continuous functions must be the limit of the sequence of conditional mathematical expectations of a sequence of random vectors.

At last, we should state the fact that, these results in this paper can be easily extended to the cases of multivariate continuous functions based on the methods in section 2.